\documentclass[11pt]{amsart}
\usepackage[mathscr]{euscript}
\usepackage{amssymb}

\usepackage{amscd}
\usepackage{amsmath, amssymb}
\usepackage{amsfonts}
\usepackage[colorlinks,linkcolor=blue,citecolor=blue, pdfstartview=FitH]{hyperref}
\usepackage[all]{xy}
\usepackage{tikz}

\numberwithin{equation}{section}

\theoremstyle{plain}
\newtheorem{thm}{Theorem}[section]
\newtheorem{theorem}[thm]{Theorem}

\newtheorem{corollary}[thm]{Corollary}
\newtheorem{proposition}[thm]{Proposition}

\theoremstyle{definition}

\newtheorem{question}[thm]{Question}

\newtheorem{remark}[thm]{Remark}

\newtheorem{definition}[thm]{Definition}

\newtheorem{defn-thm}[thm]{Definition-Theorem}

\newcommand{\im}{{ \mathrm{im}\,}}

\newcommand{\C}{{ \mathbb{C} }}
\newcommand{\Z}{{ \mathbb{Z} }}

\newcommand{\p}{{ \partial }}
\newcommand{\pb}{{ \bar{\partial} }}

\newcommand{\wg}{{\,\wedge\,}}

\begin{document}
\title{The structure of deformed double complexes on the Iwasawa manifold}
\author{Yan Hu, Wei Xia}
\address{Yan Hu, Mathematical Science Research Center, Chongqing University of Technology, Chongqing 400054, China} \email{huyan1128@126.com}
\address{Wei Xia, Mathematical Science Research Center, Chongqing University of Technology, Chongqing 400054, China} \email{xiaweiwei3@126.com, xiawei@cqut.edu.cn}

\thanks{This work is supported by the National Natural Science Foundation of China No. 11901590, the Natural Science Foundation of Chongqing (China) No. CSTB2022NSCQ-MSX0876 and Scientific Research Foundation of Chongqing University of Technology.}
\date{\today}

\begin{abstract}
The Kuranishi family of the Iwasawa manifold give rise naturally to a family of (deformed) double complexes. By using the structure theorem of double complexes due to Stelzig and Qi-Khovanov, we show there are exactly $3$ isomorphism types in this family and determine explicitly structures of these $3$ types. As an application, we computed the Fr\"olicher spectral sequence for each fiber in the Kuranishi family of the Iwasawa manifold.
\end{abstract}

\maketitle
\section{Introduction}
Given a complex manifold $X$, the space of (complex) differential forms $A^{\bullet}(X)$ admits a natural double complex structure with differential given by $d=\p+\pb$. From this double complex, denoted by
\begin{equation}\label{eq-Adppb-0}
(A^{\bullet}(X),d=\p+\pb),
\end{equation}
the Dolbeault/Bott-Chern/Aeppli cohomology and the Fr\"olicher spectral sequence can be derived. These objects has been extensively studied in recent years (see e.g. \cite{AK17a,RZ15,RYY19,AT13,Kas13,PSU21,CFG91,CFUG97}). In order to attain a deeper comprehension of these complex structure invariants, it is necessary to study the structure of $(A^{\bullet}(X),d=\p+\pb)$. In fact, according to the structure theorem of double complexes due to Stelzig \cite{Ste21} and Qi-Khovanov \cite{KQ20}, any double complex can be decomposed into a direct sum of indecomposable ones (i.e. squares and zigzags) and the Dolbeault/Bott-Chern/Aeppli cohomology, the Fr\"olicher spectral sequence can be determined by counting zigzags of specified types.

Instead of a single complex manifold $X$, we are more interested in how the double complex $(A^{\bullet}(X),d=\p+\pb)$ vary as the complex structure on $X$ deforms.

Let $\pi: (\mathcal{X}, X)\to (B,0)$ be a small deformation of the compact complex manifold $X$ such that for each $t\in B$ the complex structure on $X_t$ is represented by Beltrami differential $\phi(t)\in A^{0,1}(X,T_X^{1,0})$. For any $t\in B$, there is a double complex
\begin{equation}\label{eq-Adphit-0}
(A^{\bullet}(X),d_{\phi(t)}=\p+\pb_{\phi(t)}),
\end{equation}
where $\pb_{\phi(t)}=\pb-\mathcal{L}_{\phi(t)}^{1,0}$ and $\mathcal{L}_{\phi(t)}^{1,0}=i_{\phi(t)}\p-\p i_{\phi(t)}$. \eqref{eq-Adphit-0} may be regarded as an analytic family of deformations of the usual double complex \eqref{eq-Adppb-0} because it is reduced to \eqref{eq-Adppb-0} when $t=0$. The double complex \eqref{eq-Adphit-0} is naturally related to $(A^{\bullet}(X_t),d=\p_t+\pb_t)$: there is an isomorphism of filtered complex (see Subsection \ref{sec-2.2})
\[
e^{i_{\phi(t)}}=\sum_{k=0}^\infty i_{\phi(t)}^k:(F^\bullet A^{\bullet}(X),d_{\phi(t)} )\to (F^\bullet A^{\bullet}(X_t),d).
\]
In particular, the Fr\"olicher spectral sequence of the double complex \eqref{eq-Adphit-0} is isomorphic to the Fr\"olicher spectral sequence of $X_t$.

The Iwasawa manifold (or more generally, complex parallelizable solvmanifolds) and its Kuranishi family has often been adopted as a test for many propositions about cohomology and the Fr\"olicher spectral sequence, see e.g. \cite{Nak75,Ang13,RZ18,Mas20,Xia19dDol}. For the Kuranishi family $\{X_t\}_{t\in \mathcal{B}}$ of the Iwasawa manifold $X=X_0$, there is a natural double sub-complex
\[
(\wedge^{\bullet}, d_{\phi(t)}=\p+\pb_{\phi(t)})
\]
of \eqref{eq-Adphit-0} for any $t\in\mathcal{B}$ which is constituted of left invariant forms on $X$ and the inclusion map is an $E_1$-isomorphism in the sense of Stelzig \cite{Ste21}.

By using the structure theorem of double complexes, Flavi \cite{Fla20} recently computed the structure of the double sub-complex $(\wedge^{\bullet}(X_t), d=\p_t+\pb_t)$
of $(A^{\bullet}(X_t),d=\p_t+\pb_t)$ for some small deformations of the Iwasawa manifold. Based on this, Flavi determined the Fr\"olicher spectral sequence $E_r^{\bullet,\bullet}$ correspondingly. We are curious about the whole picture about deformation stability of degeneration of the Fr\"olicher spectral sequence \cite{WX23}. Motivated by these, in this paper, we explicitly decompose the double complex $(\wedge^{\bullet}, d_{\phi(t)}=\p+\pb_{\phi(t)})$ as indecomposable ones (i.e. squares and zigzags) for the Kuranishi family of the Iwasawa manifold. The main result of this paper may be summarized as follows:
\begin{thm}\label{thm-main-0}
Let $\{X_t\}_{t\in \mathcal{B}}$ be the Kuranishi family of the Iwasawa manifold $X=X_0$ which is defined by Beltrami differentials $\{\phi(t)\in A^{0,1}(X,T_X^{1,0})\}_{t\in \mathcal{B}}$. Denoted by $\wedge^{\bullet,\bullet}\subset A^{\bullet,\bullet}(X)$ the space of left invariant forms. Then for any $t\in\mathcal{B}$, the structure of double complex $(\wedge^{\bullet}, d_{\phi(t)}=\p+\pb_{\phi(t)})$ is described as follows:
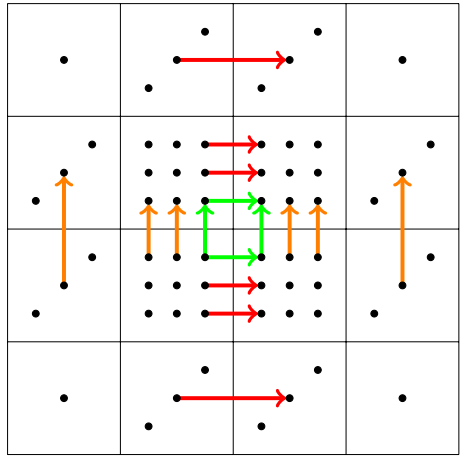
\begin{center}
\begin{tikzpicture}[scale=1.5]

\draw[help lines,black] (0,0) grid (4,4);

\draw [line width=1.5,->,red] (3/2,1/2) -- (5/2-0.03,1/2);
\draw [line width=1.5,->,orange] (1/2,3/2) -- (1/2,5/2-0.03);

\draw [line width=1.5,->,red] (7/4,5/4) -- (9/4-0.03,5/4);
\draw [line width=1.5,->,red] (7/4,3/2) -- (9/4-0.03,3/2);
\draw [line width=1.5,->,orange] (5/4,7/4) -- (5/4,9/4-0.03);
\draw [line width=1.5,->,orange] (3/2,7/4) -- (3/2,9/4-0.03);
\draw [line width=1.5,->,green] (7/4,7/4) -- (9/4-0.03,7/4);
\draw [line width=1.5,->,green] (7/4,7/4) -- (7/4,9/4-0.03);
\draw [line width=1.5,->,green] (9/4,7/4) -- (9/4,9/4-0.03);
\draw [line width=1.5,->,green] (7/4,9/4) -- (9/4-0.03,9/4);

\draw [line width=1.5,->,red] (7/4,5/2) -- (9/4-0.03,5/2);
\draw [line width=1.5,->,red] (7/4,11/4) -- (9/4-0.03,11/4);
\draw [line width=1.5,->,orange] (5/2,7/4) -- (5/2,9/4-0.03);
\draw [line width=1.5,->,orange] (11/4,7/4) -- (11/4,9/4-0.03);

\draw [line width=1.5,->,orange] (7/2,3/2) -- (7/2,5/2-0.03);
\draw [line width=1.5,->,red] (3/2,7/2) -- (5/2-0.03,7/2);

\draw [fill] (1/2,1/2) circle [radius=0.03];

\draw [fill] (5/4,1/4) circle [radius=0.03];
\draw [fill] (3/2,1/2) circle [radius=0.03];
\draw [fill] (7/4,3/4) circle [radius=0.03];

\draw [fill] (9/4,1/4) circle [radius=0.03];
\draw [fill] (5/2,1/2) circle [radius=0.03];
\draw [fill] (11/4,3/4) circle [radius=0.03];

\draw [fill] (7/2,1/2) circle [radius=0.03];

\draw [fill] (1/4,5/4) circle [radius=0.03];
\draw [fill] (1/2,3/2) circle [radius=0.03];
\draw [fill] (3/4,7/4) circle [radius=0.03];

\draw [fill] (5/4,5/4) circle [radius=0.03];
\draw [fill] (3/2,5/4) circle [radius=0.03];
\draw [fill] (7/4,5/4) circle [radius=0.03];
\draw [fill] (5/4,3/2) circle [radius=0.03];
\draw [fill] (3/2,3/2) circle [radius=0.03];
\draw [fill] (7/4,3/2) circle [radius=0.03];
\draw [fill] (5/4,7/4) circle [radius=0.03];
\draw [fill] (3/2,7/4) circle [radius=0.03];
\draw [fill] (7/4,7/4) circle [radius=0.03];

\draw [fill] (9/4,5/4) circle [radius=0.03];
\draw [fill] (5/2,5/4) circle [radius=0.03];
\draw [fill] (11/4,5/4) circle [radius=0.03];
\draw [fill] (9/4,3/2) circle [radius=0.03];
\draw [fill] (5/2,3/2) circle [radius=0.03];
\draw [fill] (11/4,3/2) circle [radius=0.03];
\draw [fill] (9/4,7/4) circle [radius=0.03];
\draw [fill] (5/2,7/4) circle [radius=0.03];
\draw [fill] (11/4,7/4) circle [radius=0.03];

\draw [fill] (13/4,5/4) circle [radius=0.03];
\draw [fill] (7/2,3/2) circle [radius=0.03];
\draw [fill] (15/4,7/4) circle [radius=0.03];

\draw [fill] (1/4,9/4) circle [radius=0.03];
\draw [fill] (1/2,5/2) circle [radius=0.03];
\draw [fill] (3/4,11/4) circle [radius=0.03];

\draw [fill] (5/4,9/4) circle [radius=0.03];
\draw [fill] (3/2,9/4) circle [radius=0.03];
\draw [fill] (7/4,9/4) circle [radius=0.03];
\draw [fill] (5/4,5/2) circle [radius=0.03];
\draw [fill] (3/2,5/2) circle [radius=0.03];
\draw [fill] (7/4,5/2) circle [radius=0.03];
\draw [fill] (5/4,11/4) circle [radius=0.03];
\draw [fill] (3/2,11/4) circle [radius=0.03];
\draw [fill] (7/4,11/4) circle [radius=0.03];

\draw [fill] (9/4,9/4) circle [radius=0.03];
\draw [fill] (5/2,9/4) circle [radius=0.03];
\draw [fill] (11/4,9/4) circle [radius=0.03];
\draw [fill] (9/4,5/2) circle [radius=0.03];
\draw [fill] (5/2,5/2) circle [radius=0.03];
\draw [fill] (11/4,5/2) circle [radius=0.03];
\draw [fill] (9/4,11/4) circle [radius=0.03];
\draw [fill] (5/2,11/4) circle [radius=0.03];
\draw [fill] (11/4,11/4) circle [radius=0.03];

\draw [fill] (13/4,9/4) circle [radius=0.03];
\draw [fill] (7/2,5/2) circle [radius=0.03];
\draw [fill] (15/4,11/4) circle [radius=0.03];

\draw [fill] (1/2,7/2) circle [radius=0.03];

\draw [fill] (5/4,13/4) circle [radius=0.03];
\draw [fill] (3/2,7/2) circle [radius=0.03];
\draw [fill] (7/4,15/4) circle [radius=0.03];

\draw [fill] (9/4,13/4) circle [radius=0.03];
\draw [fill] (5/2,7/2) circle [radius=0.03];
\draw [fill] (11/4,15/4) circle [radius=0.03];
\draw [fill] (7/2,7/2) circle [radius=0.03];

\draw [fill] (2,-1/4) node {$(i)~(t_{11},t_{12},t_{21},t_{22})= 0, D(t)=0$};
\end{tikzpicture}
\end{center}

\begin{center}
\begin{tikzpicture}[scale=1.5]
\draw[help lines,black] (0,0) grid (4,4);
\draw [line width=1.5,->,blue] (3/2,1/2) -- (5/2-0.03,1/2);
\draw [line width=1.5,->,orange] (1/2,3/2) -- (1/2,5/2-0.03);
\draw [line width=1.5,->,blue] (3/2,1/2) -- (3/2,5/4-0.03);

\draw[line width=1.5,->,purple]  (9/4,1/4) -- (9/4,5/4-0.03);
\draw [line width=1.5,->,purple] (7/4,5/4) -- (9/4-0.03,5/4);
\draw [line width=1.5,->,red] (7/4,3/2) -- (9/4-0.03,3/2);
\draw [line width=1.5,->,orange] (5/4,7/4) -- (5/4,9/4-0.03);
\draw [line width=1.5,->,orange] (3/2,7/4) -- (3/2,9/4-0.03);
\draw [line width=1.5,->,green] (7/4,7/4) -- (9/4-0.03,7/4);
\draw [line width=1.5,->,green] (7/4,7/4) -- (7/4,9/4-0.03);
\draw [line width=1.5,->,green] (9/4,7/4) -- (9/4,9/4-0.03);
\draw [line width=1.5,->,green] (7/4,9/4) -- (9/4-0.03,9/4);

\draw [line width=1.5,->,red] (7/4,5/2) -- (9/4-0.03,5/2);
\draw [line width=1.5,->,blue] (7/4,11/4) -- (9/4-0.03,11/4);
\draw [line width=1.5,->,orange] (5/2,7/4) -- (5/2,9/4-0.03);
\draw [line width=1.5,->,orange] (11/4,7/4) -- (11/4,9/4);
\draw [line width=1.5,->,blue] (7/4,11/4) -- (7/4,15/4-0.03);

\draw [line width=1.5,->,orange] (7/2,3/2) -- (7/2,5/2-0.03);
\draw [line width=1.5,->,purple] (3/2,7/2) -- (5/2-0.03,7/2);
\draw [line width=1.5,->,purple] (5/2,11/4) -- (5/2,7/2-0.03);

\draw [fill] (1/2,1/2)   circle [radius=0.03];

\draw [fill] (5/4,1/4)  circle [radius=0.03];
\draw [fill] (3/2,1/2)  circle [radius=0.03];
\draw [fill] (7/4,3/4) circle [radius=0.03];

\draw [fill] (9/4,1/4) circle [radius=0.03];
\draw [fill] (5/2,1/2) circle [radius=0.03];
\draw [fill] (11/4,3/4) circle [radius=0.03];

\draw [fill] (7/2,1/2) circle [radius=0.03];

\draw [fill] (1/4,5/4) circle [radius=0.03];
\draw [fill] (1/2,3/2) circle [radius=0.03];
\draw [fill] (3/4,7/4) circle [radius=0.03];

\draw [fill] (5/4,5/4) circle [radius=0.03];
\draw [fill] (3/2,5/4) circle [radius=0.03];
\draw [fill] (7/4,5/4) circle [radius=0.03];
\draw [fill] (5/4,3/2) circle [radius=0.03];
\draw [fill] (3/2,3/2)  circle [radius=0.03];
\draw [fill] (7/4,3/2)  circle [radius=0.03];
\draw [fill] (5/4,7/4)  circle [radius=0.03];
\draw [fill] (3/2,7/4)  circle [radius=0.03];
\draw [fill] (7/4,7/4)  circle [radius=0.03];

\draw [fill] (9/4,5/4) circle [radius=0.03];
\draw [fill] (5/2,5/4) circle [radius=0.03];
\draw [fill] (11/4,5/4) circle [radius=0.03];
\draw [fill] (9/4,3/2) circle [radius=0.03];
\draw [fill] (5/2,3/2) circle [radius=0.03];
\draw [fill] (11/4,3/2) circle [radius=0.03];
\draw [fill] (9/4,7/4) circle [radius=0.03];
\draw [fill] (5/2,7/4) circle [radius=0.03];
\draw [fill] (11/4,7/4) circle [radius=0.03];

\draw [fill] (13/4,5/4) circle [radius=0.03];;
\draw [fill] (7/2,3/2) circle [radius=0.03];
\draw [fill] (15/4,7/4) circle [radius=0.03];

\draw [fill] (1/4,9/4) circle [radius=0.03];
\draw [fill] (1/2,5/2) circle [radius=0.03];
\draw [fill] (3/4,11/4) circle [radius=0.03];

\draw [fill] (5/4,9/4) circle [radius=0.03];
\draw [fill] (3/2,9/4) circle [radius=0.03];
\draw [fill] (7/4,9/4) circle [radius=0.03];
\draw [fill] (5/4,5/2) circle [radius=0.03];
\draw [fill] (3/2,5/2) circle [radius=0.03];
\draw [fill] (7/4,5/2) circle [radius=0.03];
\draw [fill] (5/4,11/4) circle [radius=0.03];
\draw [fill] (3/2,11/4) circle [radius=0.03];
\draw [fill] (7/4,11/4) circle [radius=0.03];

\draw [fill] (9/4,9/4) circle [radius=0.03];
\draw [fill] (5/2,9/4) circle [radius=0.03];
\draw [fill] (11/4,9/4) circle [radius=0.03];
\draw [fill] (9/4,5/2) circle [radius=0.03];
\draw [fill] (5/2,5/2) circle [radius=0.03];
\draw [fill] (11/4,5/2) circle [radius=0.03];
\draw [fill] (9/4,11/4) circle [radius=0.03];
\draw [fill] (5/2,11/4) circle [radius=0.03];
\draw [fill] (11/4,11/4) circle [radius=0.03];

\draw [fill] (13/4,9/4) circle [radius=0.03];
\draw [fill] (7/2,5/2) circle [radius=0.03];
\draw [fill] (15/4,11/4) circle [radius=0.03];

\draw [fill] (1/2,7/2) circle [radius=0.03];

\draw [fill] (5/4,13/4) circle [radius=0.03];
\draw [fill] (3/2,7/2) circle [radius=0.03];
\draw [fill] (7/4,15/4) circle [radius=0.03];

\draw [fill] (9/4,13/4) circle [radius=0.03];
\draw [fill] (5/2,7/2) circle [radius=0.03];
\draw [fill] (11/4,15/4) circle [radius=0.03];

\draw [fill] (7/2,7/2) circle [radius=0.03];

\draw [fill] (2,-1/4) node {$(ii)~(t_{11},t_{12},t_{21},t_{22})\neq 0, D(t)=0$};
\end{tikzpicture}\qquad
\begin{tikzpicture}[scale=1.5]
\draw[help lines,black] (0,0) grid (4,4);
\draw [line width=1.5,->,blue] (3/2,1/2) -- (5/2-0.03,1/2);
\draw [line width=1.5,->,orange] (1/2,3/2) -- (1/2,5/2-0.03);
\draw [line width=1.5,->,blue] (3/2,1/2) -- (3/2,5/4-0.03);

\draw[line width=1.5,->,purple]  (9/4,1/4) -- (9/4,5/4-0.03);
\draw [line width=1.5,->,purple] (7/4,5/4) -- (9/4-0.03,5/4);
\draw [line width=1.5,->,purple] (7/4,3/2) -- (11/4-0.03,3/2);
\draw [line width=1.5,->,orange] (5/4,7/4) -- (5/4,9/4-0.03);
\draw [line width=1.5,->,orange] (3/2,7/4) -- (3/2,9/4-0.03);
\draw [line width=1.5,->,green] (7/4,7/4) -- (9/4-0.03,7/4);
\draw [line width=1.5,->,green] (7/4,7/4) -- (7/4,9/4-0.03);
\draw [line width=1.5,->,green] (9/4,7/4) -- (9/4,9/4-0.03);
\draw [line width=1.5,->,green] (7/4,9/4) -- (9/4-0.03,9/4);
\draw [line width=1.5,->,purple] (11/4,3/4) -- (11/4,3/2-0.03);

\draw [line width=1.5,->,blue] (5/4,5/2) -- (9/4-0.03,5/2);
\draw [line width=1.5,->,blue] (7/4,11/4) -- (9/4-0.03,11/4);
\draw [line width=1.5,->,orange] (5/2,7/4) -- (5/2,9/4-0.03);
\draw [line width=1.5,->,orange] (11/4,7/4) -- (11/4,9/4-0.03);
\draw [line width=1.5,->,blue] (7/4,11/4) -- (7/4,15/4-0.03);
\draw [line width=1.5,->,blue] (5/4,5/2) -- (5/4,13/4-0.03);

\draw [line width=1.5,->,orange] (7/2,3/2) -- (7/2,5/2-0.03);
\draw [line width=1.5,->,purple] (3/2,7/2) -- (5/2-0.03,7/2);
\draw [line width=1.5,->,purple] (5/2,11/4) -- (5/2,7/2-0.03);

\draw [fill] (1/2,1/2)   circle [radius=0.03];

\draw [fill] (5/4,1/4)  circle [radius=0.03];
\draw [fill] (3/2,1/2) circle [radius=0.03];
\draw [fill] (7/4,3/4) circle [radius=0.03];

\draw [fill] (9/4,1/4) circle [radius=0.03];
\draw [fill] (5/2,1/2) circle [radius=0.03];
\draw [fill] (11/4,3/4) circle [radius=0.03];

\draw [fill] (7/2,1/2) circle [radius=0.03];

\draw [fill] (1/4,5/4) circle [radius=0.03];
\draw [fill] (1/2,3/2) circle [radius=0.03];
\draw [fill] (3/4,7/4) circle [radius=0.03];

\draw [fill] (5/4,5/4) circle [radius=0.03];
\draw [fill] (3/2,5/4) circle [radius=0.03];
\draw [fill] (7/4,5/4) circle [radius=0.03];
\draw [fill] (5/4,3/2) circle [radius=0.03];
\draw [fill] (3/2,3/2)  circle [radius=0.03];
\draw [fill] (7/4,3/2)  circle [radius=0.03];
\draw [fill] (5/4,7/4)  circle [radius=0.03];
\draw [fill] (3/2,7/4)  circle [radius=0.03];
\draw [fill] (7/4,7/4)  circle [radius=0.03];

\draw [fill] (9/4,5/4) circle [radius=0.03];
\draw [fill] (5/2,5/4) circle [radius=0.03];
\draw [fill] (9/4,3/2) circle [radius=0.03];
\draw [fill] (5/2,3/2) circle [radius=0.03];
\draw [fill] (11/4,3/2) circle [radius=0.03];
\draw [fill] (9/4,7/4) circle [radius=0.03];
\draw [fill] (5/2,7/4) circle [radius=0.03];
\draw [fill] (11/4,7/4) circle [radius=0.03];
\draw [fill] (11/4,5/4) circle [radius=0.03];

\draw [fill] (13/4,5/4) circle [radius=0.03];
\draw [fill] (7/2,3/2) circle [radius=0.03];
\draw [fill] (15/4,7/4) circle [radius=0.03];

\draw [fill] (1/4,9/4) circle [radius=0.03];
\draw [fill] (1/2,5/2) circle [radius=0.03];
\draw [fill] (3/4,11/4) circle [radius=0.03];

\draw [fill] (5/4,9/4) circle [radius=0.03];
\draw [fill] (3/2,9/4) circle [radius=0.03];
\draw [fill] (7/4,9/4) circle [radius=0.03];
\draw [fill] (5/4,5/2) circle [radius=0.03];
\draw [fill] (3/2,5/2) circle [radius=0.03];
\draw [fill] (7/4,5/2) circle [radius=0.03];
\draw [fill] (5/4,11/4) circle [radius=0.03];
\draw [fill] (3/2,11/4) circle [radius=0.03];
\draw [fill] (7/4,11/4) circle [radius=0.03];

\draw [fill] (9/4,9/4) circle [radius=0.03];
\draw [fill] (5/2,9/4) circle [radius=0.03];
\draw [fill] (11/4,9/4) circle [radius=0.03];
\draw [fill] (9/4,5/2) circle [radius=0.03];
\draw [fill] (5/2,5/2) circle [radius=0.03];
\draw [fill] (11/4,5/2) circle [radius=0.03];
\draw [fill] (9/4,11/4) circle [radius=0.03];
\draw [fill] (5/2,11/4) circle [radius=0.03];
\draw [fill] (11/4,11/4) circle [radius=0.03];

\draw [fill] (13/4,9/4) circle [radius=0.03];
\draw [fill] (7/2,5/2) circle [radius=0.03];
\draw [fill] (15/4,11/4) circle [radius=0.03];

\draw [fill] (1/2,7/2) circle [radius=0.03];

\draw [fill] (5/4,13/4) circle [radius=0.03];
\draw [fill] (3/2,7/2) circle [radius=0.03];
\draw [fill] (7/4,15/4) circle [radius=0.03];

\draw [fill] (9/4,13/4) circle [radius=0.03];
\draw [fill] (5/2,7/2) circle [radius=0.03];
\draw [fill] (11/4,15/4) circle [radius=0.03];

\draw [fill] (7/2,7/2) circle [radius=0.03];

\draw [fill] (2,-1/4) node {$(iii)~(t_{11},t_{12},t_{21},t_{22})\neq 0, D(t)\neq0$};
\end{tikzpicture}
\end{center}
where $D(t)=t_{11}t_{22}-t_{21}t_{12}$.
\end{thm}
\begin{remark}
Compared to the result of Flavi \cite{Fla20}, we see that the structure of double complex $(\wedge^{\bullet}, d_{\phi(t)}=\p+\pb_{\phi(t)})$ seems less complicated than that of $(\wedge^{\bullet}(X_t), d=\p_t+\pb_t)$: the former has exactly $3$ isomorphism types while the latter has at least $4$ isomorphism types. According to the computation of the Bott-Chern cohomology by Angella \cite{Ang13}, there should be at least $5$ isomorphism types of $(\wedge^{\bullet}(X_t), d=\p_t+\pb_t)$.
\end{remark}

As a consequence, we can easily compute the deformed Dolbeault/Bott-Chern/Aeppli cohomology \cite{Xia19dDol,Xia19dBC,HX24} and the Fr\"olicher spectral sequence of $X_t$. In view of these (see Table 1 and Table 2 in Section \ref{sec-5}), we get the following
\begin{corollary}\label{coro-0}
The following statements hold:
\begin{enumerate}
  \item Any small deformation of the Iwasawa manifold satisfies: $E_2=E_\infty$. There are small deformations (i.e. case $(iii)$) of the Iwasawa manifold for which $E_1=E_\infty$ holds;
  \item For the Kuranishi family of the Iwasawa manifold, both $\dim\im d_1^{p,q}$ and $\dim \frac{E_1^{p,q}}{E_2^{p,q}}=\dim\im d_1^{p,q}+\dim\im d_1^{p-1,q}$ are upper semi-continuous function of $t\in\mathcal{B}$ (under the analytic Zariski topology);
  \item The varying of $\sum_{i=0}^q(-1)^{q-i}\dim E_2^{p,i}$ may be neither upper nor lower semi-continuous.
\end{enumerate}
\end{corollary}
It is known \cite{Kas15} that for complex parallelizable solvmanifolds, $E_2=E_\infty$ always holds, that is, the Fr\"olicher spectral sequence degenerates at the second page. Moreover, it is conjectured by Kasuya-Stelzig \cite{KS23} that this holds for all complex parallelisable manifold. We note that for case $(ii)$ and $(iii)$, the deformation $X_t$ of the Iwasawa manifold is not complex parallelisable since $\dim H^{2,0}_{\pb}(X_t)=2\neq3$.

A special case of a result due to Flenner \cite{Fle81} says that the alternating sum of Hodge numbers: $\sum_{i=0}^q(-1)^{q-i}\dim E_1^{p,i}$, is upper semi-continuous for any $p$ and $q$. From $(3)$ of Corollary \ref{coro-0}, we see that this statement does not hold for $\sum_{i=0}^q(-1)^{q-i}\dim E_2^{p,i}$. It is perhaps worth noting that the upper or lower semi-continuous property of $\sum_{i=0}^q(-1)^{q-i}\dim E_2^{p,i}$ for all $p,q$ implies that of $\dim E_2^{p,q}$ but not conversely, compare \cite{Mas20,Fla20}.

\begin{question}Let $\{X_t\}_{t\in B}$ be a complex analytic family and $(E_r^{\bullet,\bullet}, d_r^{\bullet,\bullet})$ the Fr\"olicher spectral sequence on $X_t$. Is $\dim\im d_1^{p,q}$ or $\dim\im \frac{E_1^{p,q}}{E_2^{p,q}}=\dim\im d_1^{p,q}+\dim\im d_1^{p-1,q}$ upper semi-continuous function of $t\in B$?
\end{question}
\section{Preliminary}
\subsection{The double complex $(A^{\bullet}(X),d=\p+\pb)$}
Let $X$ be a complex manifold and $A^{\bullet}(X)$ the space of (complex valued) smooth forms on $X$. The de Rham complex $(A^{\bullet}(X),d)$ naturally has a double complex structure given by
\[
A^{\bullet}(X)=\bigoplus_{p\in \mathbb{N}} A^{p,\bullet-p}(X),\quad d=\p+\pb,
\]
where $A^{p,\bullet-p}(X)$ is the space of $(p,\bullet-p)$-forms on $X$. This double complex $(A^{\bullet}(X),d=\p+\pb)$ induce a filtered complex
\[
(F^\bullet A^{\bullet}(X),d),~\text{with}~F^pA^{\bullet}(X)=\oplus_{\lambda\geq p} A^{\lambda,\bullet-\lambda}(X),
\]
whose associated spectral sequence $(E_r^{p,q}(X),d_r(X))$ is called the \emph{Fr\"olicher spectral sequence} of $X$.
\subsection{The deformed double complex $(A^{\bullet}(X),d_{\phi(t)}=\p+\pb_{\phi(t)})$}\label{sec-2.2}
Now, let
\[
\pi: (\mathcal{X}, X)\to (B,0)
\]
be a complex analytic family (i.e. a surjective, proper and homomorphic submersion between complex manifolds). We make the following conventions which will be used throughout this paper:
\begin{itemize}
  \item $X=\pi^{-1}(0)$ and $X_t=\pi^{-1}(t)$ for any $t\in B$;
  \item since only small deformations will be considered, we may assume $B$ is a small polydisc in some Euclidean space;
  \item for any $t\in B$, the complex structure on $X_t=\pi^{-1}(t)$ is represented by (the Beltrami differential) $\phi(t)\in A^{0,1}(X,T^{1,0})$.
\end{itemize}
By using the Beltrami differential $\phi(t)$, we can relate smooth forms on $X$ to those on $X_t$. In fact, the exponential operator \cite{Tod89,Cle05,LSY09,LRY15}
\[
e^{i_{\phi(t)}}:=\sum_{k=0}^{\infty} \frac{i_{\phi(t)}^k}{k!},
\]
where $i_{\phi(t)}$ is the contraction operator, satisfy (see \cite{FM06,FM09,WZ20})
\begin{equation}\label{eq-weidingchang}
e^{i_{\phi(t)}}~: F^pA^{k}(X)\longrightarrow F^pA^{k}(X_t).
\end{equation}
Set (see e.g. \cite{LRY15,FM06,Xia19deri,Ma05})
\begin{equation}\label{eq-LRY15}
d_{\phi(t)}:=e^{-i_{\phi(t)}}de^{i_{\phi(t)}} = \p + \pb_{\phi(t)}=\p + \pb-\mathcal{L}_{\phi(t)}^{1,0}~,
\end{equation}
where $\pb_{\phi(t)}:=\pb-\mathcal{L}_{\phi(t)}^{1,0}$ and $\mathcal{L}_{\phi(t)}^{1,0}:=i_{\phi(t)}\p-\p i_{\phi(t)}$ is the Lie derivative. Then $(A^{\bullet}(X),d_{\phi(t)}=\p+\pb_{\phi(t)})$ constitutes a family of \emph{deformed double complex}: when $t=0$, it is reduced to $(A^{\bullet}(X),d=\p+\pb)$.

There is an intimate relation between $(A^{\bullet}(X),d_{\phi(t)}=\p+\pb_{\phi(t)})$ and $(A^{\bullet}(X_t),d=\p_t+\pb_t)$: indeed, as observed in \cite{WX23},
\[
(F^\bullet A^{\bullet}(X),d_{\phi(t)} )\xrightarrow{e^{i_{\phi(t)}}} (F^\bullet A^{\bullet}(X_t),d)
\]
is an isomorphism of filtered complexes which induces the following isomorphism
\[
E_r^{\bullet,\bullet}(X,d_{\phi(t)})\cong  E_r^{\bullet,\bullet}(X_t),
\]
where $E_r^{\bullet,\bullet}(X,d_{\phi(t)})$ is the spectral sequence associated to $(F^\bullet A^{\bullet}(X),d_{\phi(t)} )$.

\subsection{The structure theorem for double complexes}
In this subsection, we briefly recall the structure theorem for double complexes. The readers is referred to \cite{Ste21,KQ20} for more information.

A \emph{double complex} (over $\C$) is a bigraded $\C$-vector space
\[
A=\oplus_{p,q\in\Z}A^{p,q}
\]
together with linear maps
\[
\p_1:A^{p,q}\longrightarrow A^{p+1,q},\quad \p_2:A^{p,q}\longrightarrow A^{p,q+1},\quad p,q\in\Z,
\]
such that $d=\p_1+\p_2$ is a differential, i.e. $(\p_1+\p_2)^2=0$. A nonzero double complex $A$ is called \emph{indecomposable} if there exist no nontrivial decomposition $A=A_1\oplus A_2$. It is possible to write down all indecomposable double complexes. To describe these, we make the following conventions:
\begin{enumerate}
  \item All drawn components are one-dimensional;
  \item All drawn maps are isomorphisms;
  \item All components and maps not drawn are zero.
\end{enumerate}
With these in minds, indecomposable double complexes may be listed as following
\begin{itemize}
  \item Squares:
\[
\xymatrix{
  A^{p-1,q}  \ar[r]^-{\p_1} & A^{p,q}  \\
  A^{p-1,q-1} \ar[u]_{\p_2}\ar[r]^{\p_1} & A^{p,q-1}\ar[u]_{\p_2}, }
\]
  \item Zigzags:
\begin{align*}
&A^{p,q}(\text{dot}),~\xymatrix{
  A^{p,q+1}   \\
  A^{p,q}\ar[u]_{\p_2} },~ \xymatrix{ A^{p,q}\ar[r]^-{\p_1} & A^{p+1,q} },~\xymatrix{
  A^{p,q+1}   \\
  A^{p,q}\ar[u]_{\p_2}\ar[r]^-{\p_1} & A^{p+1,q} },\\
  &\xymatrix{ A^{p-1,q}\ar[r]^-{\p_1} & A^{p,q}\\
  A^{p-1,q-1}\ar[u]_{\p_2} },~\xymatrix{ A^{p-1,q}\ar[r]^-{\p_1} & A^{p,q} &\\
  & A^{p,q-1}\ar[u]_{\p_2}\ar[r]^-{\p_1}  & A^{p+1,q-1} },\cdots
\end{align*}
\end{itemize}
For a square or a zigzag $A$, the \emph{shape} is defined to be
\[
S(A)=\{(p,q)\in\Z^2\mid A^{p,q}\neq0\}.
\]
For any bounded $S\subset \Z^2$, there is an (essentially unique) double complex $A$ with
\[
A^{p,q}=
\left\{
\begin{array}{rcl}
0, & (p,q)\in S, \\[5pt]
\C, & (p,q)\notin S,
\end{array}
\right.
\]
and $(\p_1)^{p,q}=(-1)^q\text{Id}, (\p_2)^{p,q}=\text{Id}$. This double complex will be denoted by $C(S)$. Now we are ready to state the structure theorem of double complexes due to Stelzig \cite[Thm.\,3.4]{Ste21} and Qi-Khovanov \cite{KQ20}:
\begin{theorem}[Structure theorem of double complexes]\label{thm-str-thm-doublecomplex}
For every bounded double complex $A$, there exist unique cardinal numbers $\text{mult}_S(A)$ and a (non-unique) isomorphism
\[
A\cong\oplus_{S}C(S)^{\oplus\text{mult}_S(A)},
\]
where $S\subset \Z^2$ runs over the set of all shapes of squares and zigzags.
\end{theorem}
We will apply this theorem to the double complex $(\wedge^{\bullet}, d_{\phi(t)}=\p+\pb_{\phi(t)})$ for the Kuranishi family of the Iwasawa manifold. In this case, each component $\wedge^{p,q}$ of the double complex is of finite dimensional so that the decomposition process terminates at finite steps.

Since the (Dolbeault, Bott-Chern, Aeppli) cohomology groups and the Fr\"olicher spectral sequence are most easy to compute for squares and zigzags, as a consequence of Theorem \ref{thm-str-thm-doublecomplex}, we have the following (see \cite{RZ20} for some related discussion)
\begin{theorem}[Proposition $6$ and Lemma $8$ in \cite{Ste21}]\label{thm-cohomology-Er}
Let $(A=\bigoplus A^{p,q},d=\p_1+\p_2)$ be a bounded double complex with a decomposition into elementary complexes with pairwise distinct support $\varphi:\bigoplus A_i \xrightarrow{\cong} A$.
\begin{enumerate}
  \item
	For every $(p,q)\in\Z^2$, the maps induced by $\varphi$
	\[
	\bigoplus_{\substack{A_i\text{ zigzag}\\(p,q)\in S(A_i)\\(p-1,q),(p+1,q)\not\in S(A_i)}}H_{\p_1}^{p,q}(A_i)\longrightarrow H^{p,q}_{\p_1}(A)
	\]
	\[
	\bigoplus_{\substack{A_i\text{ zigzag}\\(p,q)\in S(A_i)\\(p,q-1),(p,q+1)\not\in S(A_i)}}H_{\p_2}^{p,q}(A_i)\longrightarrow H^{p,q}_{\p_2}(A)	
	\]
	are isomorphisms.
	\item
	For every $(p,q)\in\Z^2$, the maps induced by $\varphi$
	\[
	\bigoplus_{\substack{A_i\text{ zigzag}\\(p,q)\in S(A_i)\\(p+1,q),(p,q+1)\not\in S(A_i)}}H_{BC}^{p,q}(A_i)\longrightarrow H^{p,q}_{BC}(A)
	\]
	\[
	\bigoplus_{\substack{A_i\text{ zigzag}\\(p,q)\in S(A_i)\\(p-1,q),(p,q-1)\not\in S(A_i)}}H_{A}^{p,q}(A_i)\longrightarrow H^{p,q}_A(A)	
	\]
	are isomorphisms.
  \item For any $i=1,2$, let $({}_iE_r^{p,q},{}_id_r^{p,q})$ be the spectral sequence associated to the filtered complex
      \[
      (F_i^\bullet A^k,d),\quad~\text{with}~F_1^\bullet A^k=\bigoplus_{p\geq \bullet} A^{p,k-p},~F_2^\bullet A^k=\bigoplus_{q\geq \bullet} A^{k-q,q}.
      \]
      Then there is an equality
		\[
		\textit{mult}_{S_{i,r}^{p,q}}(A)=\dim\im{}_id_r^{p,q},
		\]
where $\textit{mult}_{S_{i,r}^{p,q}}(A)$ is the number of zigzags (with length $2r$) in $\bigoplus_i A_i$ which has one and only one outgoing arrow at $(p,q)$, $i$ is $1$ (if this arrow is $\p_1$) or $2$ (if this arrow is $\p_2$).
\end{enumerate}
\end{theorem}
\begin{definition}
Let $r\in\mathbb{N}$. A morphism $f:A\to B$ of double
complexes is called an \emph{$E_r$-isomorphism} \cite[Def.\,10]{Ste21} if it induces an isomorphism on the
$r$-th page of both Fr\"olicher spectral sequences (where the $0$-th page is defined to be the complex itself).
\end{definition}

\section{The Kuranishi family of the Iwasawa manifold and the double complex $(\wedge^{\bullet}, d_{\phi(t)}=\p+\pb_{\phi(t)})$ }
\subsection{Iwasawa manifold and its Kuranishi family}
Let $G$ be the complex Lie group defined by
\[
G := \left\{
\left(
\begin{array}{ccc}
 1 & z^1 & z^3 \\
 0 &  1  & z^2 \\
 0 &  0  &  1
\end{array}
\right) \in \mathrm{GL}(3;\mathbb{C}) \mid z^1,\,z^2,\,z^3 \in\mathbb{C} \right\}\cong \mathbb{C}^3,
\]
where the product is the one induced by matrix multiplication. This is usually called the \emph{Heisenberg group}. Consider the discrete subgroup $\Gamma$ defined by
\[
\Gamma := \left\{
\left(
\begin{array}{ccc}
 1 & \omega^1 & \omega^3 \\
 0 &  1  & \omega^2 \\
 0 &  0  &  1
\end{array}
\right) \in G \mid \omega^1,\,\omega^2,\,\omega^3 \in\mathbb{Z}[\sqrt{-1}] \right\},
\]
The quotient $X=\Gamma\setminus G$ is called the \emph{Iwasawa manifold}. A basis of $H^0(X,\Omega^1)$ is given by
\[
\varphi^1 = d z^1,~ \varphi^2 = d z^2,~ \varphi^3 = d z^3-z^1\,d z^2,
\]
and a dual basis $\theta^1, \theta^2, \theta^3\in H^0(X,T_X^{1,0})$ is given by
\[
\theta^1=\frac{\partial}{\partial z^1},~\theta^2=\frac{\partial}{\partial z^2} + z^1\frac{\partial}{\partial z^3},~\theta^3=\frac{\partial}{\partial z^3}.
\]
Notice that $\varphi^1, \varphi^2, \varphi^3$ are left invariant forms and $\theta^1, \theta^2, \theta^3$ are left invariant vector fields on $G$. Denote by $(\wedge^{\bullet}, d=\p+\pb)$ the double sub-complex of $(A^{\bullet}(X),d=\p+\pb)$ consisting of left invariant forms, then
\[
\wedge^{\bullet}=\oplus_{p\geq0}\wedge^{p,\bullet-p}\quad\text{and}\quad \wedge^{1,0}=\C\{\varphi^1, \varphi^2, \varphi^3\},
\]
where $\wedge^{p,\bullet-p}$ is the vector space of left invariant $(p,\bullet-p)$-forms. The structure equation in terms of the basis $\{\varphi^1,\varphi^2,\varphi^3\}$ is given by
\begin{equation}\label{eq-str-eq}
d\varphi^1=0,\quad d\varphi^2=0,\quad d\varphi^3=-\varphi^{12}.
\end{equation}
$X$ is equipped with the Hermitian metric $\sum_{i=1}^3\varphi^i\otimes\bar{\varphi}^i$. The Beltrami differential of the Kuranishi family of $X$ is
\[
\phi(t) = \sum_{i=1}^3\sum_{\lambda=1}^2t_{i\lambda}\theta^i\bar{\varphi}^{\lambda} - D(t)\theta^3\bar{\varphi}^{3},~\text{with}~D(t)=t_{11}t_{22}-t_{21}t_{12},
\]
and the Kuranishi space of $X$ is
\[
\mathcal{B}=\{t=(t_{11}, t_{12}, t_{21}, t_{22}, t_{31}, t_{32})\in \mathbb{C}^6\mid |t_{i\lambda}|<\epsilon, i=1, 2, 3, \lambda=1,2 \},
\]
where $\epsilon>0$ is sufficiently small. Since $\phi(t)$ is also left invariant, we see that for any $t\in\mathcal{B}$,
\[
(\wedge^{\bullet}, d_{\phi(t)}=\p+\pb_{\phi(t)})\stackrel{\iota}\hookrightarrow(A^{\bullet}(X),d_{\phi(t)}=\p+\pb_{\phi(t)}),
\]
is a double sub-complex. It is known that the Dolbeault cohomology of $(A^{\bullet}(X),d_{\phi(t)}=\p+\pb_{\phi(t)})$ can be computed by using only left forms \cite{Xia20g}. As a consequence, $\iota$ is an $E_1$-isomorphism for each $t\in\mathcal{B}$. To determine the decomposition of $(\wedge^{\bullet}, d=\p+\pb_{\phi(t)})$ indicated in Theorem \ref{thm-str-thm-doublecomplex}, we need first to compute the differentials $\p$ and $\pb_{\phi(t)}$ so as to find all possible shapes of squares and zigzags. The former is easy to get by using the structure equation \eqref{eq-str-eq}. So it is left to compute $\pb_{\phi(t)}$.

\subsection{Compute $\pb_{\phi(t)}$ on the subcomplex of invariant forms}\label{subsec-compute-pbphit}
Set
\[
\phi_1=\sum_{i=1}^3\sum_{\lambda=1}^2t_{i\lambda}\theta^i\bar{\varphi}^{\lambda},~ \phi_2 = -D(t)\theta^3\bar{\varphi}^{3},
\]
then by direct computations, we have
\begin{align*}
&\mathcal{L}_{\phi_1}^{1,0}\varphi^1 = \mathcal{L}_{\phi_1}^{1,0}\varphi^2=\mathcal{L}_{\phi_1}^{1,0}\bar{\varphi}^1=\mathcal{L}_{\phi_1}^{1,0}\bar{\varphi}^{2}=\mathcal{L}_{\phi_1}^{1,0}\bar{\varphi}^{3}=0,\\[5pt]
&\mathcal{L}_{\phi_1}^{1,0}\varphi^3 = \sum_{\lambda=1}^2(t_{1\lambda}\varphi^2-t_{2\lambda}\varphi^1)\wedge\bar{\varphi}^{\lambda},~\mathcal{L}_{\phi_2}^{1,0}\varphi^i =\mathcal{L}_{\phi_2}^{1,0}\bar{\varphi}^{i}=0,~i=1,2,3.
\end{align*}
In particular,
\begin{equation}\label{eq-plie}
\pb_{\phi(t)}=\pb-\mathcal{L}_{\phi_1}^{1,0}.
\end{equation}

\subsubsection{$\pb_{\phi(t)}:\wedge^{1,0}\to\wedge^{1,1}=\C\{\varphi^{i\bar{j}}\mid 1\leq i,j\leq3\}$}\label{subsubsec-10to11}
Recall that $\wedge^{1,0}=\C\{\varphi^1,\varphi^2,\varphi^3 \}$ and $\wedge^{1,1}=\C\{\varphi^{i\bar{j}}\mid 1\leq i,j\leq3\}$.
\begin{align*}
\pb_{\phi(t)}\varphi^1&=(\pb-\mathcal{L}_{\phi_1}^{1,0})\varphi^1=0,\\
\pb_{\phi(t)}\varphi^2&=(\pb-\mathcal{L}_{\phi_1}^{1,0})\varphi^2=0,\\
\pb_{\phi(t)}\varphi^3&=(\pb-\mathcal{L}_{\phi_1}^{1,0})\varphi^3
=-\mathcal{L}_{\phi_1}^{1,0}\varphi^3\\
&=-\sum_{\lambda=1}^2(t_{1\lambda}\varphi^2-t_{2\lambda}\varphi^1)\wg\varphi^{\bar{\lambda}}\\
&=t_{21}\varphi^{1\bar{1}}-t_{11}\varphi^{2\bar{1}}+t_{22}\varphi^{1\bar{2}}-t_{12}\varphi^{2\bar{2}}.
\end{align*}
We see that
\[
\dim\im\pb_{\phi(t)}\cap\wedge^{1,1}=\left\{
\begin{array}{rcl}
0, &(t_{11},t_{12},t_{21},t_{22}) =0 \\[5pt]
1, &(t_{11},t_{12},t_{21},t_{22})\neq 0
\end{array}
\right. .
\]
\subsubsection{$\pb_{\phi(t)}:\wedge^{0,1}\to\wedge^{0,2}$}\label{subsubsec-01to02}
Recall that $\wedge^{0,1}=\C\{\varphi^{\bar{1}},\varphi^{\bar{2}},
\varphi^{\bar{3}}\}$ and $\wedge^{0,2}=\C\{\varphi^{\overline{12}}
,\varphi^{\overline{13}},\varphi^{\overline{23}}\}$.
\begin{align*}
\pb_{\phi(t)}\varphi^{\bar{1}}&=(\pb-\mathcal{L}_{\phi_1}^{1,0})\varphi^{\bar{1}}=0,\\
\pb_{\phi(t)}\varphi^{\bar{2}}&=(\pb-\mathcal{L}_{\phi_1}^{1,0})\varphi^{\bar{2}}=0,\\
\pb_{\phi(t)}\varphi^{\bar{3}}&
=(\pb-\mathcal{L}_{\phi_1}^{1,0})\varphi^{\bar{3}}=-\varphi^{\overline{12}}\neq 0.
\end{align*}
We see that $\dim\im\pb_{\phi(t)}\cap\wedge^{0,2}=1$.
\subsubsection{$\pb_{\phi(t)}:\wedge^{2,0}\to\wedge^{2,1}$}\label{subsubsec-20to21}
Recall that $\wedge^{2,0}=\C\{\varphi^{12},\varphi^{13},
\varphi^{23}\}$ and $\wedge^{2,1}=\C\{\varphi^{ij\overline{k}}\mid 1\leq i,j,k\leq3, i<j\}$.
\begin{align*}
\pb_{\phi(t)}\varphi^{12}&=(\pb-\mathcal{L}_{\phi_1}^{1,0})\varphi^{12}=0,\\
\pb_{\phi(t)}\varphi^{13}&=(\pb-\mathcal{L}_{\phi_1}^{1,0})\varphi^{13}
=t_{11}\varphi^{12\bar{1}}+t_{12}\varphi^{12\bar{2}}\\
\pb_{\phi(t)}\varphi^{23}&=(\pb-\mathcal{L}_{\phi_1}^{1,0})\varphi^{23}
=t_{21}\varphi^{12\bar{1}}+t_{22}\varphi^{12\bar{2}}.
\end{align*}
Let $\sigma_0=a_1\varphi^{13}+a_2\varphi^{23}$, then
\begin{align*}
\pb_{\phi(t)}\sigma_0=&a_1(t_{11}\varphi^{12\bar{1}}+t_{12}\varphi^{12\bar{2}})+a_2(t_{21}\varphi^{12\bar{1}}+t_{22}\varphi^{12\bar{2}})\\
=&(a_1t_{11}+a_2t_{21})\varphi^{12\bar{1}}+(a_1t_{12}+a_2t_{22})\varphi^{12\bar{2}}.
\end{align*}
Set
\[
T=\left(\begin{array}{cc}
    t_{11} & t_{21} \\
    t_{12} & t_{22}
\end{array}\right),
\]
then
\begin{enumerate}
  \item When $(t_{11},t_{12},t_{21},t_{22})=0$, we have $\text{rank}(T)=0$ and $\dim\im\pb_{\phi(t)}\cap\wedge^{2,1}=0$.
  \item When $(t_{11},t_{12},t_{21},t_{22})\neq 0,D(t)=0$, we have $\text{rank}(T)=1$ and $\dim\im\pb_{\phi(t)}\cap\wedge^{2,1}=1$.
  \item When $(t_{11},t_{12},t_{21},t_{22})\neq 0,D(t)\neq 0$, we have $\text{rank}(T)=2$ and $\dim\im\pb_{\phi(t)}\cap\wedge^{2,1}=2$.
\end{enumerate}
\subsubsection{$\pb_{\phi(t)}:\wedge^{1,1}\to\wedge^{1,2}$}\label{subsubsec-11to12}
Recall that $\wedge^{1,1}=\C\{\varphi^{1\bar{1}},\varphi^{1\bar{2}},\varphi^{1\bar{3}},\varphi^{2\bar{1}},\varphi^{2\bar{2}},\varphi^{2\bar{3}},\varphi^{3\bar{1}},
\varphi^{3\bar{2}},\varphi^{3\bar{3}}\}$ and $\wedge^{1,2}=\C\{\varphi^{i\overline{jk}}\mid 1\leq i,j,k\leq3, j<k\}$.
\begin{align*}
\pb_{\phi(t)}\varphi^{1\bar{1}}&=(\pb-\mathcal{L}_{\phi_1}^{1,0})\varphi^{1\bar{1}}=0,\\
\pb_{\phi(t)}\varphi^{1\bar{2}}&=(\pb-\mathcal{L}_{\phi_1}^{1,0})\varphi^{1\bar{2}}=0,\\
\pb_{\phi(t)}\varphi^{2\bar{1}}&=(\pb-\mathcal{L}_{\phi_1}^{1,0})\varphi^{2\bar{1}}=0,\\
\pb_{\phi(t)}\varphi^{2\bar{2}}&=(\pb-\mathcal{L}_{\phi_1}^{1,0})\varphi^{2\bar{2}}=0,\\
\pb_{\phi(t)}\varphi^{1\bar{3}}&=(\pb-\mathcal{L}_{\phi_1}^{1,0})\varphi^{1\bar{3}}=\varphi^{1\overline{12}},\\ \pb_{\phi(t)}\varphi^{2\bar{3}}&=(\pb-\mathcal{L}_{\phi_1}^{1,0})\varphi^{2\bar{3}}=\varphi^{2\overline{12}},\\
\pb_{\phi(t)}\varphi^{3\bar{1}}&=(\pb-\mathcal{L}_{\phi_1}^{1,0})\varphi^{3\bar{1}}=-\mathcal{L}_{\phi_1}^{1,0}\varphi^{3\bar{1}}\\
&=-\sum_{\lambda=1}^2(t_{1\lambda}\varphi^2-t_{2\lambda}\varphi^1)\wg\varphi^{\bar{\lambda}}\wg \varphi^{\bar{1}}\\
&=-t_{22}\varphi^{1\overline{12}}+t_{12}\varphi^{2\overline{12}},\\
\pb_{\phi(t)}\varphi^{3\bar{2}}&=(\pb-\mathcal{L}_{\phi_1}^{1,0})\varphi^{3\bar{2}}=-\mathcal{L}_{\phi_1}^{1,0}\varphi^{3\bar{2}}\\
&=-\sum_{\lambda=1}^2(t_{1\lambda}\varphi^2-t_{2\lambda}\varphi^1)\wg\varphi^{\bar{\lambda}}\wg \varphi^{\bar{2}}\\
&=t_{21}\varphi^{1\overline{12}}-t_{11}\varphi^{2\overline{12}}.
\end{align*}
Let $\sigma_0=a_1\varphi^{1\bar{3}}+a_2\varphi^{2\bar{3}}+a_3\varphi^{3\bar{1}}+a_4\varphi^{3\bar{2}}$, then
\begin{align*}
\pb_{\phi(t)}\sigma_0&=a_1\varphi^{1\overline{12}}+a_2\varphi^{2\overline{12}}+a_3(-t_{22}\varphi^{1\overline{12}}+t_{12}\varphi^{2\overline{12}})
+a_4(t_{21}\varphi^{1\overline{12}}-t_{11}\varphi^{2\overline{12}})\\
&=(a_1-a_3t_{22}+a_4t_{21})\varphi^{1\overline{12}}+(a_2+a_3t_{12}-a_4t_{11})\varphi^{2\overline{12}}.
\end{align*}
Set
\[
T=\left(\begin{array}{cccc}
   1 & 0 &  -t_{22} & t_{21} \\
    0 & 1 & t_{12} & -t_{11}
\end{array}\right),
\]
we see that the rank of
\[
\pb_{\phi(t)}:\C\{\varphi^{1\bar{3}}, \varphi^{2\bar{3}}, \varphi^{3\bar{1}}, \varphi^{3\bar{2}}\}\longrightarrow\C\{\varphi^{1\overline{12}}, \varphi^{2\overline{12}}\},
\]
is equal to $\text{rank}(T)=2$. Similarly, from
\begin{align*}
\pb_{\phi(t)}\varphi^{3\bar{3}}&=(\pb-\mathcal{L}_{\phi_1}^{1,0})\varphi^{3\bar{3}}=-\varphi^3\wg\pb\varphi^{\bar{3}}
-\mathcal{L}_{\phi_1}^{1,0}\varphi^3\wg \varphi^{\bar{3}}\\
&=\varphi^{3\overline{12}}-\sum_{\lambda=1}^2(t_{1\lambda}\varphi^2-t_{2\lambda}\varphi^1)\wg\varphi^{\bar{\lambda}}\wg \varphi^{\bar{3}}\\
&=\varphi^{3\overline{12}}-t_{11}\varphi^{2\overline{13}}+t_{21}\varphi^{1\overline{13}}-t_{12}\varphi^{2\overline{23}}+t_{22}\varphi^{1\overline{23}},
\end{align*}
we see that there the rank of
\[
\pb_{\phi(t)}:\C\{\varphi^{3\bar{3}}\}\longrightarrow\C\{\varphi^{1\overline{13}}, \varphi^{1\overline{23}}, \varphi^{2\overline{13}}, \varphi^{2\overline{23}}, \varphi^{3\overline{12}}\},
\]
is equal to $1$.
\subsubsection{$\pb_{\phi(t)}:\wedge^{0,2}\to\wedge^{0,3}$}\label{subsubsec-02to03}
Recall that $\wedge^{0,2}=\C\{\varphi^{\overline{12}}, \varphi^{\overline{13}}, \varphi^{\overline{23}}\}$ and $\wedge^{0,3}=\C\{\varphi^{\overline{123}}\}$.
\begin{align*}
 \pb_{\phi(t)}\varphi^{\overline{12}}&=(\pb-\mathcal{L}_{\phi_1}^{1,0})\varphi^{\overline{12}}=0,\\
  \pb_{\phi(t)}\varphi^{\overline{13}}&=(\pb-\mathcal{L}_{\phi_1}^{1,0})\varphi^{\overline{13}}=0,\\
\pb_{\phi(t)}\varphi^{\overline{23}}&=(\pb-\mathcal{L}_{\phi_1}^{1,0})\varphi^{\overline{23}}=0.
 \end{align*}
We see that $\dim\im\pb_{\phi(t)}\cap\wedge^{0,3}=0$.
\subsubsection{$\pb_{\phi(t)}:\wedge^{3,0}\to\wedge^{3,1}$}\label{subsubsec-30to31}
Recall that $\wedge^{3,0}=\C\{\varphi^{123}\}$ and $\wedge^{3,1}=\C\{\varphi^{123\overline{1}},\varphi^{123\overline{2}},\varphi^{123\overline{3}}\}$.
\[
\pb_{\phi(t)}\varphi^{123}=(\pb-\mathcal{L}_{\phi_1}^{1,0})\varphi^{123}=-\mathcal{L}_{\phi_1}^{1,0}\varphi^{123}=0.
\]
We see that $\dim\im\pb_{\phi(t)}\cap\wedge^{3,1}=0$.
\subsubsection{$\pb_{\phi(t)}:\wedge^{2,1}\to\wedge^{2,2}$}\label{subsubsec-21to22}
Recall that $\wedge^{2,1}=\C\{\varphi^{ij\bar{k}}\mid 1\leq i,j,k\leq3, i<j\}$ and $\wedge^{2,2}=\C\{\varphi^{ij\overline{kl}}\mid 1\leq i,j,k,l\leq3, k<l\}$.
\begin{align*}
\pb_{\phi(t)}\varphi^{12\bar{1}}&=(\pb-\mathcal{L}_{\phi_1}^{1,0})\varphi^{12\bar{1}}=0,\\
\pb_{\phi(t)}\varphi^{12\bar{2}}&=(\pb-\mathcal{L}_{\phi_1}^{1,0})\varphi^{12\bar{2}}=0,\\
\pb_{\phi(t)}\varphi^{12\bar{3}}&=(\pb-\mathcal{L}_{\phi_1}^{1,0})\varphi^{12\bar{3}}=-\varphi^{12\overline{12}},
\end{align*}
\begin{align*}
\pb_{\phi(t)}\varphi^{13\bar{1}}&=(\pb-\mathcal{L}_{\phi_1}^{1,0})\varphi^{13\bar{1}}
=-\mathcal{L}_{\phi_1}^{1,0}\varphi^{13\bar{1}}
=-t_{12}\varphi^{12\overline{12}},\\
\pb_{\phi(t)}\varphi^{13\bar{2}}&=(\pb-\mathcal{L}_{\phi_1}^{1,0})\varphi^{13\bar{2}}
=-\mathcal{L}_{\phi_1}^{1,0}\varphi^{13\bar{2}}
=t_{11}\varphi^{12\overline{12}},\\
\pb_{\phi(t)}\varphi^{23\bar{1}}&=(\pb-\mathcal{L}_{\phi_1}^{1,0})\varphi^{23\bar{1}}
=-\mathcal{L}_{\phi_1}^{1,0}\varphi^{23\bar{1}}
=-t_{22}\varphi^{12\overline{12}},\\
\pb_{\phi(t)}\varphi^{23\bar{2}}&=(\pb-\mathcal{L}_{\phi_1}^{1,0})\varphi^{23\bar{2}}
=-\mathcal{L}_{\phi_1}^{1,0}\varphi^{23\bar{2}}
=t_{21}\varphi^{12\overline{12}}.
\end{align*}
Let $\sigma_0=a_1\varphi^{12\bar{3}}+a_2\varphi^{13\bar{1}}+a_3\varphi^{13\bar{2}}+a_4\varphi^{23\bar{1}}+a_5\varphi^{23\bar{2}}$, then
\[
\pb_{\phi(t)}\sigma_0=(-a_1-a_2t_{12}+a_3t_{11}-a_4t_{22}+a_5t_{21})\varphi^{12\overline{12}}.
\]
Set
\[
T=\left(\begin{array}{ccccc}
   -1 & -t_{21} & t_{11}  & -t_{22} & t_{21}
   \end{array}\right).
\]
We see that the rank of
\[
\pb_{\phi(t)}:\C\{\varphi^{12\bar{3}}, \varphi^{13\bar{1}}, \varphi^{13\bar{2}}, \varphi^{23\bar{1}},\varphi^{23\bar{2}}\}\longrightarrow\C\{\varphi^{12\overline{12}}\}
\]
is equal to $\text{rank}(T)=1$. Similarly, from
\begin{align*}
\pb_{\phi(t)}\varphi^{13\bar{3}}&=(\pb-\mathcal{L}_{\phi_1}^{1,0})\varphi^{13\bar{3}}\\
&=-\varphi^{13\overline{12}}+\varphi^1\wg[(t_{11}\varphi^2-t_{21}\varphi^1)\wg\varphi^{\overline{13}}+(t_{12}\varphi^2-t_{22}\varphi^1)\wg\varphi^{\overline{23}}]\\
&=-\varphi^{13\overline{12}}+t_{11}\varphi^{12\overline{13}}+t_{12}\varphi^{12\overline{23}},
\end{align*}
we see that the rank of
\[
\pb_{\phi(t)}:\C\{\varphi^{13\bar{3}}\}\longrightarrow\C\{\varphi^{13\overline{12}}, \varphi^{12\overline{13}}, \varphi^{12\overline{23}}\},
\]
is equal to $1$.

From
\begin{align*}
\pb_{\phi(t)}\varphi^{23\bar{3}}&=(\pb-\mathcal{L}_{\phi_1}^{1,0})\varphi^{23\bar{3}}\\
&=-\varphi^{23\overline{12}}+\varphi^2\wg[\sum_{\lambda=1}^2(t_{1\lambda}\varphi^2-t_{2\lambda}\varphi^1)\wg\varphi^{\bar{\lambda}}]\wg \varphi^{\bar{3}}\\
&=-\varphi^{23\overline{12}}+t_{21}\varphi^{12\overline{13}}+t_{22}\varphi^{12\overline{23}},
\end{align*}
the rank of
\[
\pb_{\phi(t)}:\C\{\varphi^{23\bar{3}}\}\longrightarrow\C\{\varphi^{23\overline{12}}, \varphi^{12\overline{13}}, \varphi^{12\overline{23}}\},
\]
is equal to $1$.
\subsubsection{$\pb_{\phi(t)}:\wedge^{1,2}\to\wedge^{1,3}$}\label{subsubsec-12to13}
Recall that $\wedge^{1,2}=\C\{\varphi^{i\overline{jk}}\mid 1\leq i,j,k\leq3, j<k\}$ and $\wedge^{1,3}=\C\{\varphi^{1\overline{123}},\varphi^{2\overline{123}},\varphi^{3\overline{123}},\}$.
\begin{align*}
\pb_{\phi(t)}\varphi^{1\overline{12}}&=(\pb-\mathcal{L}_{\phi_1}^{1,0})\varphi^{1\overline{12}}=0,\\
\pb_{\phi(t)}\varphi^{1\overline{13}}&=(\pb-\mathcal{L}_{\phi_1}^{1,0})\varphi^{1\overline{13}}=0,\\
\pb_{\phi(t)}\varphi^{1\overline{23}}&=(\pb-\mathcal{L}_{\phi_1}^{1,0})\varphi^{1\overline{23}}=0,\\
\pb_{\phi(t)}\varphi^{2\overline{12}}&=(\pb-\mathcal{L}_{\phi_1}^{1,0})\varphi^{2\overline{12}}=0,\\
\pb_{\phi(t)}\varphi^{2\overline{13}}&=(\pb-\mathcal{L}_{\phi_1}^{1,0})\varphi^{2\overline{13}}=0,\\
\pb_{\phi(t)}\varphi^{2\overline{23}}&=(\pb-\mathcal{L}_{\phi_1}^{1,0})\varphi^{2\overline{23}}=0,
\end{align*}
\begin{align*}
\pb_{\phi(t)}\varphi^{3\overline{12}}&=(\pb-\mathcal{L}_{\phi_1}^{1,0})\varphi^{3\overline{12}}=-\mathcal{L}_{\phi_1}^{1,0}\varphi^{3\overline{12}}\\
&=-\sum_{\lambda=1}^2(t_{1\lambda}\varphi^2-t_{2\lambda}\varphi^1)\wg\varphi^{\bar{\lambda}}\wg\varphi^{\overline{12}}=0,\\
\pb_{\phi(t)}\varphi^{3\overline{13}}&=(\pb-\mathcal{L}_{\phi_1}^{1,0})\varphi^{3\overline{13}}\\
&=-\sum_{\lambda=1}^2(t_{1\lambda}\varphi^2-t_{2\lambda}\varphi^1)\wg\varphi^{\bar{\lambda}}\wg\varphi^{\overline{13}}\\
&=t_{12}\varphi^{2\overline{123}}-t_{22}\varphi^{1\overline{123}},\\
\pb_{\phi(t)}\varphi^{3\overline{23}}&=(\pb-\mathcal{L}_{\phi_1}^{1,0})\varphi^{3\overline{23}}\\
&=-\sum_{\lambda=1}^2(t_{1\lambda}\varphi^2-t_{2\lambda}\varphi^1)\wg\varphi^{\bar{\lambda}}\wg\varphi^{\overline{23}}\\
&=t_{21}\varphi^{1\overline{123}}-t_{11}\varphi^{2\overline{123}}.
\end{align*}
Let $\sigma_0=a_1\varphi^{3\overline{13}}+a_2\varphi^{3\overline{23}}$, then
\begin{align*}
\pb_{\phi(t)}\sigma_0&=a_1(t_{12}\varphi^{2\overline{123}}-t_{22}\varphi^{1\overline{123}})+a_2(t_{21}\varphi^{1\overline{123}}-t_{11}\varphi^{2\overline{123}})\\
&=(-a_1t_{22}-a_2t_{21})\varphi^{1\overline{123}}+(a_1t_{12}-a_2t_{11})\varphi^{2\overline{123}}.
\end{align*}
Set
\[
T=\left(\begin{array}{cc}
    -t_{22} & t_{21} \\
    t_{12} & -t_{11}
\end{array}\right),
\]
then
\begin{enumerate}
  \item When $(t_{11},t_{12},t_{21},t_{22})=0$, we have $\text{rank}(T)=0$ and the rank of $\pb_{\phi(t)}:\C\{\varphi^{3\overline{13}}, \varphi^{3\overline{23}}\}\longrightarrow\C\{\varphi^{1\overline{123}}, \varphi^{2\overline{123}}\}$ is $0$.
  \item When $(t_{11},t_{12},t_{21},t_{22})\neq 0,D(t)=0$, we have $\text{rank}(T)=1$, and the rank of $\pb_{\phi(t)}:\C\{\varphi^{3\overline{13}}, \varphi^{3\overline{23}}\}\longrightarrow\C\{\varphi^{1\overline{123}}, \varphi^{2\overline{123}}\}$ is $1$.
  \item When $(t_{11},t_{12},t_{21},t_{22})\neq 0,D(t)\neq 0$, we have $\text{rank}(T)=2$, and the rank of $\pb_{\phi(t)}:\C\{\varphi^{3\overline{13}}, \varphi^{3\overline{23}}\}\longrightarrow\C\{\varphi^{1\overline{123}}, \varphi^{2\overline{123}}\}$ is $2$.
\end{enumerate}

\subsubsection{$\pb_{\phi(t)}:\wedge^{3,1}\to\wedge^{3,2}$}\label{subsubsec-31to32}
Recall that $\wedge^{3,1}=\C\{\varphi^{123\overline{1}},\varphi^{123\overline{2}},\varphi^{123\overline{3}},\}$ and $\wedge^{3,2}=\C\{\varphi^{123\overline{12}},\varphi^{123\overline{13}},\varphi^{123\overline{23}}\}$.
\begin{align*}
\pb_{\phi(t)}\varphi^{123\overline{1}}&=(\pb-\mathcal{L}_{\phi_1}^{1,0})\varphi^{123\overline{1}}=-\mathcal{L}_{\phi_1}^{1,0}\varphi^{123\overline{1}}\\
&=-\varphi^{12}\wg[\sum_{\lambda=1}^2(t_{1\lambda}\varphi^2-t_{2\lambda}\varphi^1)\wg\varphi^{\bar{\lambda}}]\wg\varphi^{\bar{1}}=0,
\end{align*}
\begin{align*}
\pb_{\phi(t)}\varphi^{123\overline{2}}&=(\pb-\mathcal{L}_{\phi_1}^{1,0})\varphi^{123\overline{2}}=-\mathcal{L}_{\phi_1}^{1,0}\varphi^{123\overline{2}}\\
&=-\varphi^{12}\wg[\sum_{\lambda=1}^2(t_{1\lambda}\varphi^2-t_{2\lambda}\varphi^1)\wg\varphi^{\bar{\lambda}}]\wg\varphi^{\bar{2}}=0,\\
\pb_{\phi(t)}\varphi^{123\overline{3}}&=(\pb-\mathcal{L}_{\phi_1}^{1,0})\varphi^{123\overline{3}}\\
&=\varphi^{123\overline{12}}-\varphi^{12}\wg[\sum_{\lambda=1}^2(t_{1\lambda}\varphi^2-t_{2\lambda}\varphi^1)\wg\varphi^{\bar{\lambda}}]\wg\varphi^{\bar{3}}\\
&=\varphi^{123\overline{12}}.
\end{align*}
We see that the rank of $\pb_{\phi(t)}:\C\{\varphi^{123\bar{3}}\}\longrightarrow\C\{\varphi^{123\overline{12}}\}$ is $1$.
\subsubsection{$\pb_{\phi(t)}:\wedge^{2,2}\to\wedge^{2,3}$}\label{subsubsec-22to23}
Recall that $\wedge^{2,2}=\C\{\varphi^{ij\overline{kl}}\mid 1\leq i,j,k\leq3, i<j,k<l\}$ and $\wedge^{2,3}=\C\{\varphi^{12\overline{123}},\varphi^{13\overline{123}},\varphi^{23\overline{123}}\}$.
\begin{align*}
\pb_{\phi(t)}\varphi^{12\overline{12}}&=(\pb-\mathcal{L}_{\phi_1}^{1,0})\varphi^{12\overline{12}}=0,\\
\pb_{\phi(t)}\varphi^{12\overline{13}}&=(\pb-\mathcal{L}_{\phi_1}^{1,0})\varphi^{12\overline{13}}=0,\\
\pb_{\phi(t)}\varphi^{12\overline{23}}&=(\pb-\mathcal{L}_{\phi_1}^{1,0})\varphi^{12\overline{23}}=0,\\
\pb_{\phi(t)}\varphi^{13\overline{12}}&=(\pb-\mathcal{L}_{\phi_1}^{1,0})\varphi^{13\overline{12}}=-\mathcal{L}_{\phi_1}^{1,0}\varphi^{13\overline{12}}\\
&=\varphi^1\wg[\sum_{\lambda=1}^2(t_{1\lambda}\varphi^2-t_{2\lambda}\varphi^1)\wg\varphi^{\bar{\lambda}}]\wg\varphi^{\overline{12}}=0,\\
\pb_{\phi(t)}\varphi^{23\overline{12}}&=(\pb-\mathcal{L}_{\phi_1}^{1,0})\varphi^{23\overline{12}}=-\mathcal{L}_{\phi_1}^{1,0}\varphi^{23\overline{12}}\\
&=\varphi^2\wg[\sum_{\lambda=1}^2(t_{1\lambda}\varphi^2-t_{2\lambda}\varphi^1)\wg\varphi^{\bar{\lambda}}]\wg\varphi^{\overline{12}}=0,
\end{align*}
\begin{align*}
\pb_{\phi(t)}\varphi^{13\overline{13}}&=(\pb-\mathcal{L}_{\phi_1}^{1,0})\varphi^{13\overline{13}}\\
&=\varphi^1\wg[\sum_{\lambda=1}^2(t_{1\lambda}\varphi^2-t_{2\lambda}\varphi^1\wg\varphi^{\bar{\lambda}}]\wg\varphi^{\overline{13}}\\
&=-t_{12}\varphi^{12\overline{123}},\\
\pb_{\phi(t)}\varphi^{13\overline{23}}&=(\pb-\mathcal{L}_{\phi_1}^{1,0})\varphi^{13\overline{23}}\\
&=\varphi^1\wg[\sum_{\lambda=1}^2(t_{1\lambda}\varphi^2-t_{2\lambda}\varphi^1\wg\varphi^{\bar{\lambda}}]\wg\varphi^{\overline{23}}\\
&=t_{11}\varphi^{12\overline{123}},
\end{align*}
\begin{align*}
\pb_{\phi(t)}\varphi^{23\overline{13}}&=(\pb-\mathcal{L}_{\phi_1}^{1,0})\varphi^{23\overline{13}}=-\mathcal{L}_{\phi_1}^{1,0}\varphi^{23\overline{13}}\\
&=\varphi^2\wg[\sum_{\lambda=1}^2(t_{1\lambda}\varphi^2-t_{2\lambda}\varphi^1)\wg\varphi^{\bar{\lambda}}]\wg\varphi^{\overline{13}}\\
&=-t_{22}\varphi^{12\overline{123}},\\
\pb_{\phi(t)}\varphi^{23\overline{23}}&=(\pb-\mathcal{L}_{\phi_1}^{1,0})\varphi^{23\overline{23}}=-\mathcal{L}_{\phi_1}^{1,0}\varphi^{23\overline{23}}\\
&=\varphi^2\wg[\sum_{\lambda=1}^2(t_{1\lambda}\varphi^2-t_{2\lambda}\varphi^1)\wg\varphi^{\bar{\lambda}}]\wg\varphi^{\overline{23}}\\
&=t_{21}\varphi^{12\overline{123}}.
\end{align*}
Let $\sigma_0=a_1\varphi^{13\overline{13}}+a_2\varphi^{13\overline{23}}+a_3\varphi^{23\overline{13}}+a_4\varphi^{23\overline{23}}$, then
\begin{align*}
 \pb_{\phi(t)}\sigma_0&=-a_1t_{12}\varphi^{12\overline{123}}+a_2t_{11}\varphi^{12\overline{123}}-a_3t_{22}\varphi^{12\overline{123}}+a_4a_1t_{21}\varphi^{12\overline{123}}\\
 &=(-a_1t_{12}+a_2t_{11}-a_3t_{22}+a_4t_{21})\varphi^{12\overline{123}}.
\end{align*}
Set
\[
T=\left(\begin{array}{cccc}
   -t_{12} & t_{11} & -t_{22} & t_{21}
\end{array}\right),
\]
then
\begin{enumerate}
  \item When $(t_{11},t_{12},t_{21},t_{22})=0$, we have $\text{rank}(T)=0$ and the rank of $ \pb_{\phi(t)}:\C\{\varphi^{13\overline{13}}, \varphi^{13\overline{23}},\varphi^{23\overline{13}}, \varphi^{23\overline{23}}\}\longrightarrow\C\{\varphi^{12\overline{123}}\}$ is $0$.
  \item When $(t_{11},t_{12},t_{21},t_{22})\neq 0,D(t)=0$, we have $\text{rank}(T)=1$ and the rank of $ \pb_{\phi(t)}:\C\{\varphi^{13\overline{13}}, \varphi^{13\overline{23}},\varphi^{23\overline{13}}, \varphi^{23\overline{23}}\}\longrightarrow\C\{\varphi^{12\overline{123}}\}$ is $1$.
  \item When $(t_{11},t_{12},t_{21},t_{22})\neq 0,D(t)\neq 0$, we have $\text{rank}(T)=1$ and the rank of $ \pb_{\phi(t)}:\C\{\varphi^{13\overline{13}}, \varphi^{13\overline{23}},\varphi^{23\overline{13}}, \varphi^{23\overline{23}}\}\longrightarrow\C\{\varphi^{12\overline{123}}\}$ is $1$.
\end{enumerate}

\subsubsection{$\pb_{\phi(t)}:\wedge^{3,2}\to\wedge^{3,3}$}\label{subsubsec-32to33}
Recall that $\wedge^{3,2}=\C\{\varphi^{123\overline{12}},\varphi^{123\overline{13}},\varphi^{123\overline{23}}\}$ and $\wedge^{3,3}=\C\{\varphi^{123\overline{123}}\}$.
\begin{align*}
\pb_{\phi(t)}\varphi^{123\overline{12}}&=(\pb-\mathcal{L}_{\phi_1}^{1,0})\varphi^{123\overline{12}}=-\mathcal{L}_{\phi_1}^{1,0}\varphi^{123\overline{12}}\\
&=-\varphi^{12}\wg[\sum_{\lambda=1}^2(t_{1\lambda}\varphi^2-t_{2\lambda}\varphi^1)\wg\varphi^{\bar{\lambda}}]\wg\varphi^{\overline{12}}=0,\\
\pb_{\phi(t)}\varphi^{123\overline{13}}&=(\pb-\mathcal{L}_{\phi_1}^{1,0})\varphi^{123\overline{13}}\\
&=-\varphi^{12}\wg[\sum_{\lambda=1}^2(t_{1\lambda}\varphi^2-t_{2\lambda}\varphi^1)\wg\varphi^{\bar{\lambda}}]\wg\varphi^{\overline{13}}=0,\\
\pb_{\phi(t)}\varphi^{123\overline{23}}&=(\pb-\mathcal{L}_{\phi_1}^{1,0})\varphi^{123\overline{23}}\\
&=-\varphi^{12}\wg[\sum_{\lambda=1}^2(t_{1\lambda}\varphi^2-t_{2\lambda}\varphi^1)\wg\varphi^{\bar{\lambda}}]\wg\varphi^{\overline{23}}=0.
\end{align*}
We see that $\dim\im\pb_{\phi(t)}\cap\wedge^{3,3}=0$.

\section{The decomposition of $(\wedge^{\bullet}, d_{\phi(t)}=\p+\pb_{\phi(t)})$}
According to Theorem \ref{thm-str-thm-doublecomplex}, we may decompose the double complex $(\wedge^{\bullet}, d_{\phi(t)}=\p+\pb_{\phi(t)})$ as direct sums of squares and zigzags. To realize this decomposition, we need to choose suitable basis of each components $\wedge^{p,q}$ so that they are compatible with the actions of differentials $\p, \pb_{\phi(t)}$. This means that if $\varphi_1,\cdots, \varphi_k$ is chosen as a basis for $\wedge^{p,q}$, then $\pb_{\phi(t)}\varphi_1,\cdots, \pb_{\phi(t)}\varphi_k$ should be contained in the basis we have chosen for $\wedge^{p,q+1}$ and $\p\varphi_1,\cdots, \p\varphi_k$ should also be contained in the basis we have chosen for $\wedge^{p+1,q}$. For $t=0$, this has already done by Angella in \cite{Ang13} and the result is pictured as follows:
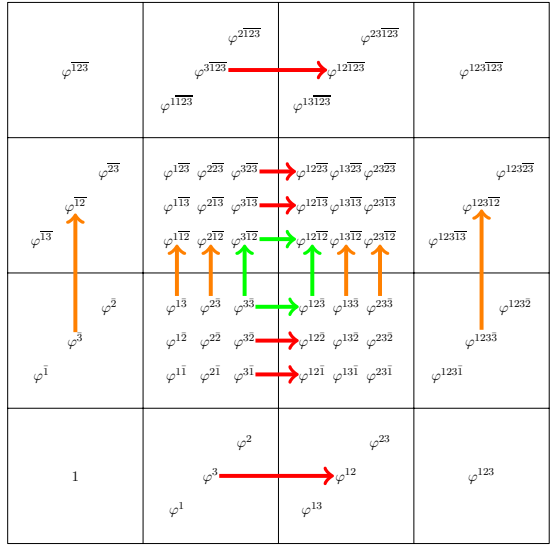
\begin{center}
\begin{tikzpicture}[scale=1.8]

\draw[help lines,black] (0,0) grid (4,4);

\draw [line width=1.5,->,red] (3/2+1/16,1/2) -- (5/2-0.09,1/2);
\draw [line width=1.5,->,orange] (1/2,3/2+1/16) -- (1/2,5/2-0.06);

\draw [line width=1.5,->,red] (7/4+0.08,5/4) -- (9/4-0.1,5/4);
\draw [line width=1.5,->,red] (7/4+0.08,3/2) -- (9/4-0.1,3/2);
\draw [line width=1.5,->,orange] (5/4,7/4+0.08) -- (5/4,9/4-0.04);
\draw [line width=1.5,->,orange] (3/2,7/4+0.08) -- (3/2,9/4-0.04);
\draw [line width=1.5,->,green] (7/4+0.08,7/4) -- (9/4-0.1,7/4);
\draw [line width=1.5,->,green] (7/4,7/4+.08) -- (7/4,9/4-0.04);
\draw [line width=1.5,->,green] (9/4,7/4+0.08) -- (9/4,9/4-0.04);
\draw [line width=1.5,->,green] (7/4+0.11,9/4) -- (9/4-0.12,9/4);

\draw [line width=1.5,->,red] (7/4+0.11,5/2) -- (9/4-0.12,5/2);
\draw [line width=1.5,->,red] (7/4+0.11,11/4) -- (9/4-0.12,11/4);
\draw [line width=1.5,->,orange] (5/2,7/4+0.08) -- (5/2,9/4-0.04);
\draw [line width=1.5,->,orange] (11/4,7/4+0.08) -- (11/4,9/4-0.04);

\draw [line width=1.5,->,orange] (7/2,3/2+0.08) -- (7/2,5/2-0.03);
\draw [line width=1.5,->,red] (3/2+0.13,7/2) -- (5/2-0.14,7/2);

\draw [fill] (1/2,1/2)   node[scale=0.5] {1};

\draw [fill=red] (5/4,1/4)  node[scale=0.5] {$\varphi^1$};
\draw [fill=red] (3/2,1/2) node[scale=0.5] {$\varphi^3$};
\draw [fill=red] (7/4,3/4) node[scale=0.5] {$\varphi^2$};

\draw [fill=green] (9/4,1/4) node[scale=0.5] {$\varphi^{13}$};
\draw [fill=red] (5/2,1/2) node[scale=0.5] {$\varphi^{12}$};
\draw [fill=green] (11/4,3/4) node[scale=0.5] {$\varphi^{23}$};

\draw [fill=blue] (7/2,1/2) node[scale=0.5] {$\varphi^{123}$};

\draw [fill=red] (1/4,5/4) node[scale=0.5] {$\varphi^{\bar{1}}$};
\draw [fill=red] (1/2,3/2) node[scale=0.5] {$\varphi^{\bar{3}}$};
\draw [fill=red] (3/4,7/4) node[scale=0.5] {$\varphi^{\bar{2}}$};

\draw [fill=green] (5/4,5/4) node[scale=0.5] {$\varphi^{1\bar{1}}$};
\draw [fill=green] (3/2,5/4) node[scale=0.5] {$\varphi^{2\bar{1}}$};
\draw [fill=green] (7/4,5/4) node[scale=0.5] {$\varphi^{3\bar{1}}$};
\draw [fill=green] (5/4,3/2) node[scale=0.5] {$\varphi^{1\bar{2}}$};
\draw [fill=green] (3/2,3/2) node[scale=0.5] {$\varphi^{2\bar{2}}$};
\draw [fill=green] (7/4,3/2) node[scale=0.5] {$\varphi^{3\bar{2}}$};
\draw [fill=green] (5/4,7/4) node[scale=0.5] {$\varphi^{1\bar{3}}$};
\draw [fill=green] (3/2,7/4) node[scale=0.5] {$\varphi^{2\bar{3}}$};
\draw [fill=green] (7/4,7/4) node[scale=0.5] {$\varphi^{3\bar{3}}$};

\draw [fill=green] (9/4,5/4) node[scale=0.5] {$\varphi^{12\bar{1}}$};
\draw [fill=blue] (5/2,5/4) node[scale=0.5] {$\varphi^{13\bar{1}}$};
\draw [fill=blue] (11/4,5/4) node[scale=0.5] {$\varphi^{23\bar{1}}$};
\draw [fill=green] (9/4,3/2) node[scale=0.5] {$\varphi^{12\bar{2}}$};
\draw [fill=blue] (5/2,3/2) node[scale=0.5] {$\varphi^{13\bar{2}}$};
\draw [fill=blue] (11/4,3/2) node[scale=0.5] {$\varphi^{23\bar{2}}$};
\draw [fill=green] (9/4,7/4) node[scale=0.5] {$\varphi^{12\bar{3}}$};
\draw [fill=blue] (5/2,7/4) node[scale=0.5] {$\varphi^{13\bar{3}}$};
\draw [fill=blue] (11/4,7/4) node[scale=0.5] {$\varphi^{23\bar{3}}$};

\draw [fill=orange] (13/4,5/4) node[scale=0.5] {$\varphi^{123\bar{1}}$};
\draw [fill=orange] (7/2,3/2) node[scale=0.5] {$\varphi^{123\bar{3}}$};
\draw [fill=orange] (15/4,7/4) node[scale=0.5] {$\varphi^{123\bar{2}}$};

\draw [fill=green] (1/4,9/4) node[scale=0.5] {$\varphi^{\overline{13}}$};
\draw [fill=red] (1/2,5/2) node[scale=0.5] {$\varphi^{\overline{12}}$};
\draw [fill=green] (3/4,11/4) node[scale=0.5] {$\varphi^{\overline{23}}$};

\draw [fill=green] (5/4,9/4) node[scale=0.5] {$\varphi^{1\overline{12}}$};
\draw [fill=green] (3/2,9/4) node[scale=0.5] {$\varphi^{2\overline{12}}$};
\draw [fill=green] (7/4,9/4) node[scale=0.5] {$\varphi^{3\overline{12}}$};
\draw [fill=blue] (5/4,5/2) node[scale=0.5] {$\varphi^{1\overline{13}}$};
\draw [fill=blue] (3/2,5/2) node[scale=0.5] {$\varphi^{2\overline{13}}$};
\draw [fill=blue] (7/4,5/2) node[scale=0.5] {$\varphi^{3\overline{13}}$};
\draw [fill=blue] (5/4,11/4) node[scale=0.5] {$\varphi^{1\overline{23}}$};
\draw [fill=blue] (3/2,11/4) node[scale=0.5] {$\varphi^{2\overline{23}}$};
\draw [fill=blue] (7/4,11/4) node[scale=0.5] {$\varphi^{3\overline{23}}$};

\draw [fill=green] (9/4,9/4) node[scale=0.5] {$\varphi^{12\overline{12}}$};
\draw [fill=blue] (5/2,9/4) node[scale=0.5] {$\varphi^{13\overline{12}}$};
\draw [fill=blue] (11/4,9/4) node[scale=0.5] {$\varphi^{23\overline{12}}$};
\draw [fill=blue] (9/4,5/2) node[scale=0.5] {$\varphi^{12\overline{13}}$};
\draw [fill=orange] (5/2,5/2) node[scale=0.5] {$\varphi^{13\overline{13}}$};
\draw [fill=orange] (11/4,5/2) node[scale=0.5] {$\varphi^{23\overline{13}}$};
\draw [fill=blue] (9/4,11/4) node[scale=0.5] {$\varphi^{12\overline{23}}$};
\draw [fill=orange] (5/2,11/4) node[scale=0.5] {$\varphi^{13\overline{23}}$};
\draw [fill=orange] (11/4,11/4) node[scale=0.5] {$\varphi^{23\overline{23}}$};

\draw [fill=magenta] (13/4,9/4) node[scale=0.5] {$\varphi^{123\overline{13}}$};
\draw [fill=orange] (7/2,5/2)  node[scale=0.5] {$\varphi^{123\overline{12}}$};
\draw [fill=magenta] (15/4,11/4) node[scale=0.5] {$\varphi^{123\overline{23}}$};

\draw [fill=blue] (1/2,7/2) node[scale=0.5] {$\varphi^{\overline{123}}$};

\draw [fill=orange] (5/4,13/4) node[scale=0.5] {$\varphi^{1\overline{123}}$};
\draw [fill=orange] (3/2,7/2) node[scale=0.5] {$\varphi^{3\overline{123}}$};
\draw [fill=orange] (7/4,15/4) node[scale=0.5] {$\varphi^{2\overline{123}}$};

\draw [fill=magenta] (9/4,13/4) node[scale=0.5] {$\varphi^{13\overline{123}}$};
\draw [fill=orange] (5/2,7/2) node[scale=0.5] {$\varphi^{12\overline{123}}$};
\draw [fill=magenta] (11/4,15/4) node[scale=0.5] {$\varphi^{23\overline{123}}$};

\draw [fill] (7/2,7/2) node[scale=0.5] {$\varphi^{123\overline{123}}$};
\end{tikzpicture}
\end{center}
where signs in several places are omitted (e.g. $\varphi^{12}$) because we are only interested in the double complex structure. The same convention applies in latter subsections so long as this does not cause any trouble.
\subsection{The case that $(t_{11},t_{12},t_{21},t_{22})\neq 0, D(t)=0$}\label{subsec-(ii)}
\begin{proposition}\label{prop-D(t)=0}
There exist the following basis
\begin{itemize}
  \item $\wedge^{1,0}=\C\{ \underbrace{\varphi^{1}, \varphi^{2}}_{\ker\pb_{\phi(t)}}, \varphi^{3} \}$;
  \item $\wedge^{2,0}=\C\{\underbrace{\xi_1^{2,0}, \varphi^{13}}_{\ker\pb_{\phi(t)}}, \xi_2^{2,0}\}$;
  \item $\wedge^{1,1}=\C\{\underbrace{\pb_{\phi(t)}\varphi^{3}, \xi_1^{1,1}, \xi_2^{1,1}, \xi_3^{1,1}, \xi_4^{1,1}, \xi_5^{1,1}}_{\ker\pb_{\phi(t)}}, \varphi^{1\bar{3}}, \varphi^{2\bar{3}}, \varphi^{3\bar{3}} \}$;
  \item $\wedge^{2,1}=\C\{\underbrace{\p\xi_4^{1,1}, \p\xi_5^{1,1}, \xi_1^{2,1}, \xi_2^{2,1}, \xi_3^{2,1}, \xi_4^{2,1}}_{\ker\pb_{\phi(t)}}, \varphi^{12\bar{3}}, \varphi^{13\bar{3}}, \varphi^{23\bar{3}} \}$;
  \item $\wedge^{1,2}=\C\{\underbrace{\pb_{\phi(t)}\varphi^{3\bar{3}}, \varphi^{1\overline{12}}, \varphi^{2\overline{12}}, \xi_1^{1,2}, \xi_2^{1,2}, \xi_3^{1,2}, \xi_4^{1,2}, \xi_5^{1,2}}_{\ker\pb_{\phi(t)}}, \xi_6^{1,2}\}$;
  \item $\wedge^{2,2}=\C\{\underbrace{\pb_{\phi(t)}\varphi^{13\bar{3}}, \pb_{\phi(t)}\varphi^{23\bar{3}}, \varphi^{12\overline{12}}, \p\xi_5^{1,2}, \p\xi_6^{1,2}, \xi_1^{2,2}, \xi_2^{2,2}, \xi_3^{2,2}}_{\ker\pb_{\phi(t)}}, \xi_4^{2,2} \}$;
  \item $\wedge^{1,3}=\C\{\pb_{\phi(t)}\xi_6^{1,2}, \xi_1^{1,3}, \varphi^{3\overline{123}} \}$,
\end{itemize}
such that
\begin{itemize}
  \item $\mathbb{C}\{\xi_1^{2,0},\xi_2^{2,0}\}=\mathbb{C}\{\varphi^{13},\varphi^{23}\}$;
  \item $\mathbb{C}\{\xi_1^{1,1},\xi_2^{1,1},\xi_3^{1,1},\pb_{\phi(t)}\varphi^3\}
      =\mathbb{C}\{\varphi^{1\bar{1}},\varphi^{1\bar{2}},
      \varphi^{2\bar{1}},\varphi^{2\bar{2}}\}$ and
  \[
   \mathbb{C}\{\xi_4^{1,1},\xi_5^{1,1}\}=
   \mathbb{C}\{\varphi^{3\bar{1}},\varphi^{3\bar{2}}\},\quad \p\xi_4^{1,1}=\pb_{\phi(t)}\xi_2^{2,0};
  \]
  \item $\mathbb{C}\{\varphi^{12\bar{3}},\xi_1^{2,1},\xi_2^{2,1},\xi_3^{2,1},
      \xi_4^{2,1}\}=\mathbb{C}\{\varphi^{12\bar{3}},\varphi^{13\bar{1}},
      \varphi^{13\bar{2}},\varphi^{23\bar{1}},\varphi^{23\bar{2}}\}$;
  \item $\mathbb{C}\{\xi_1^{1,2},\xi_2^{1,2},\xi_3^{1,2},\xi_4^{1,2},
      \pb_{\phi(t)}\varphi^{3\overline{3}}\}=
      \mathbb{C}\{\varphi^{1\overline{13}},\varphi^{1\overline{23}},
      \varphi^{2\overline{13}},\varphi^{2\overline{23}},\varphi^{3\overline{12}}\}$
      and
      \[
\mathbb{C}\{\xi_5^{1,2}, \xi_6^{1,2}\}=\C\{\varphi^{3\overline{13}}, \varphi^{3\overline{23}}\},\quad \text{with}~\pb_{\phi(t)}\xi_6^{1,2}\neq0;
\]
  \item $\mathbb{C}\{\p\xi_5^{1,2},\p\xi_6^{1,2},\pb_{\phi(t)}\varphi^{13\overline{3}},
      \pb_{\phi(t)}\varphi^{23\overline{3}}\}=\mathbb{C}\{\varphi^{12\overline{13}},
      \varphi^{12\overline{23}},\varphi^{13\overline{12}},\varphi^{23\overline{12}}\}$;
  \item $\mathbb{C}\{\xi_1^{1,3},\pb_{\phi(t)}\xi_6^{1,2}\}=
      \mathbb{C}\{\varphi^{1\overline{123}},\varphi^{2\overline{123}}\}$.
\end{itemize}
\end{proposition}
\begin{proof}From Subsection \ref{subsubsec-10to11}, we know that $\ker\pb_{\phi(t)}\cap\wedge^{1,0}=\C\{\varphi^{1}, \varphi^{2}\}$ and the rank of $\pb_{\phi(t)}:\C\{\varphi^3\}\to
\C\{\varphi^{1\bar{1}},\varphi^{2\bar{1}},\varphi^{1\bar{2}},\varphi^{2\bar{2}}\}$ is $1$.

From Subsection \ref{subsubsec-20to21}, we know that $\dim\ker\pb_{\phi(t)}\cap\wedge^{2,0}=2$ and $\pb_{\phi(t)}\varphi^{12}=0$. On the other hand, the rank of $\pb_{\phi(t)}:\C\{\varphi^{13}, \varphi^{23}\}\to
\C\{\varphi^{12\bar{1}},\varphi^{12\bar{2}}\}$ is $1$. So we may choose $\xi_1^{2,0},\xi_2^{2,0}\in\wedge^{2,0}$ such that $\mathbb{C}\{\xi_1^{2,0},\xi_2^{2,0}\}=\mathbb{C}\{\varphi^{13},\varphi^{23}\}$ and
\[
\pb_{\phi(t)}\xi_1^{2,0}=0,\quad \pb_{\phi(t)}\xi_2^{2,0}\neq 0.
\]
Assume $\xi_2^{2,0}=a\varphi^{13}+b\varphi^{23}$, then
\begin{align*}
\pb_{\phi(t)}\xi_2^{2,0}&=(at_{11}+bt_{21})\varphi^{12\bar{1}}
+(at_{12}+bt_{22})\varphi^{12\bar{2}}\\
&=\p[-(at_{11}+bt_{21})\varphi^{3\bar{1}}
-(at_{12}+bt_{22})\varphi^{3\bar{2}}].
\end{align*}
Set $\xi_4^{1,1}=-(at_{11}+bt_{21})\varphi^{3\bar{1}}
-(at_{12}+bt_{22})\varphi^{3\bar{2}}$, then $\p\xi_4^{1,1}=\pb_{\phi(t)}\xi_2^{2,0}$. We may choose $\xi_5^{1,1}\in\wedge^{1,1}$ so that $\mathbb{C}\{\xi_4^{1,1},\xi_5^{1,1}\}=\C\{\varphi^{3\bar{1}},\varphi^{3\bar{2}}\}$. In particular, $\mathbb{C}\{\p\xi_4^{1,1},\p\xi_5^{1,1}\}=\C\{\varphi^{12\bar{1}},\varphi^{12\bar{2}}\}$.

From Subsection \ref{subsubsec-11to12}, we know that $\dim\ker\pb_{\phi(t)}\cap\wedge^{1,1}=6$ and
\[
\ker\pb_{\phi(t)}=\C\{\varphi^{1\bar{1}},\varphi^{1\bar{2}},
\varphi^{2\bar{1}},\varphi^{2\bar{2}},\xi_4^{1,1},\xi_5^{1,1}\}.
\]
On the other hand, $\pb_{\phi(t)}\varphi^{1\bar{3}}=\varphi^{1\overline{12}}, \pb_{\phi(t)}\varphi^{2\bar{3}}=\varphi^{2\overline{12}}$ and the rank of
$\pb_{\phi(t)}:\C\{\varphi^{3\bar{3}}\}\longrightarrow\C\{\varphi^{1\overline{13}}, \varphi^{1\overline{23}}, \varphi^{2\overline{13}}, \varphi^{2\overline{23}}, \varphi^{3\overline{12}}\}$ is $1$. Since $0\neq \pb_{\phi(t)}\varphi^3\in
\C\{\varphi^{1\bar{1}},\varphi^{1\bar{2}},\varphi^{2\bar{1}},\varphi^{2\bar{2}}\}$, we may choose $\xi_1^{1,1},\xi_2^{1,1},\xi_3^{1,1}$ such that $\mathbb{C}\{\xi_1^{1,1},\xi_2^{1,1},\xi_3^{1,1},\pb_{\phi(t)}\varphi^3\}
=\mathbb{C}\{\varphi^{1\bar{1}},\varphi^{1\bar{2}},\varphi^{2\bar{1}},\varphi^{2\bar{2}}\}$.
Notice that
\[
\p\pb_{\phi(t)}\varphi^{3\bar{3}}=-\varphi^{12\overline{12}}
=-\pb_{\phi(t)}\varphi^{12\bar{3}}=\pb_{\phi(t)}\p\varphi^{3\bar{3}},
\]

From Subsection \ref{subsubsec-21to22}, we know that $\mathbb{C}\{\p\xi_4^{1,1},\p\xi_5^{1,1}\}=\C\{\varphi^{12\bar{1}},\varphi^{12\bar{2}}\}
\subset\ker\pb_{\phi(t)}$, the rank of $\pb_{\phi(t)}:\C\{\varphi^{12\bar{3}}, \varphi^{13\bar{1}}, \varphi^{13\bar{2}}, \varphi^{23\bar{1}},\varphi^{23\bar{2}}\}\longrightarrow\C\{\varphi^{12\overline{12}}\}$ is $1$, $\pb_{\phi(t)}\varphi^{13\bar{3}}$ and $\pb_{\phi(t)}\varphi^{23\bar{3}}$ are linearly independent. We may choose $\xi_1^{2,1},\xi_2^{2,1},\xi_3^{2,1},\xi_4^{2,1}\in\ker\pb_{\phi(t)}$ so that
\[
\mathbb{C}\{\varphi^{12\bar{3}},\xi_1^{2,1},\xi_2^{2,1},\xi_3^{2,1},\xi_4^{2,1}\}
=\mathbb{C}\{\varphi^{12\bar{3}},\varphi^{13\bar{1}},
\varphi^{13\bar{2}},\varphi^{23\bar{1}},\varphi^{23\bar{2}}\}.
\]

For $\wedge^{1,2}$, we already see that $\varphi^{1\overline{12}}=\pb_{\phi(t)}\varphi^{1\bar{3}}, \varphi^{2\overline{12}}=\pb_{\phi(t)}\varphi^{2\bar{3}}$ and the rank of
$
\pb_{\phi(t)}:\C\{\varphi^{3\bar{3}}\}\longrightarrow\C\{\varphi^{1\overline{13}}, \varphi^{1\overline{23}}, \varphi^{2\overline{13}}, \varphi^{2\overline{23}}, \varphi^{3\overline{12}}\}
$ is $1$. From Subsection \ref{subsubsec-12to13}, we know that $\varphi^{1\overline{12}}, \varphi^{1\overline{13}}, \varphi^{1\overline{23}}, \varphi^{2\overline{12}}, \varphi^{2\overline{13}}, \varphi^{2\overline{23}}, \varphi^{3\overline{12}}$ are in $\ker\pb_{\phi(t)}$ and the rank of $\pb_{\phi(t)}:\C\{\varphi^{3\overline{13}}, \varphi^{3\overline{23}}\}\longrightarrow\C\{\varphi^{1\overline{123}}, \varphi^{2\overline{123}}\}$ is $1$. So we may choose $\xi_1^{1,2},\cdots,\xi_6^{1,2}$ so that
$\mathbb{C}\{\xi_1^{1,2},\xi_2^{1,2},\xi_3^{1,2},\xi_4^{1,2},
\pb_{\phi(t)}\varphi^{3\overline{3}}\}=
\mathbb{C}\{\varphi^{1\overline{13}},\varphi^{1\overline{23}},
\varphi^{2\overline{13}},\varphi^{2\overline{23}},\varphi^{3\overline{12}}\}$
and
\[
\mathbb{C}\{\xi_5^{1,2}, \xi_6^{1,2}\}=\C\{\varphi^{3\overline{13}}, \varphi^{3\overline{23}}\},~\text{with}~
\pb_{\phi(t)}\xi_5^{1,2}=0,~\pb_{\phi(t)}\xi_6^{1,2}\neq0.
\]

For $\wedge^{2,2}$, we already see that
\[
\pb_{\phi(t)}\varphi^{13\bar{3}},\pb_{\phi(t)}\varphi^{23\bar{3}}\in\C
\{\varphi^{13\overline{12}},\varphi^{12\overline{13}},\varphi^{12\overline{23}},
\varphi^{23\overline{12}}\},
\]
are linearly independent. From Subsection \ref{subsubsec-22to23}, we know that $\varphi^{12\overline{12}}, \varphi^{12\overline{13}}, \varphi^{12\overline{23}}, \varphi^{13\overline{12}}, \varphi^{23\overline{12}}\in\ker\pb_{\phi(t)}$ and the rank of
$\pb_{\phi(t)}:\C\{\varphi^{13\overline{13}}, \varphi^{13\overline{23}},\varphi^{23\overline{13}}, \varphi^{23\overline{23}}\}\longrightarrow\C\{\varphi^{12\overline{123}}\}$ is $1$. Notice that $\C\{\varphi^{12\overline{13}},\varphi^{12\overline{23}} \}=
\mathbb{C}\{\p\xi_5^{1,2}, \p\xi_6^{1,2}\}$ and
\[
\C\{\pb_{\phi(t)}\varphi^{13\bar{3}},\pb_{\phi(t)}\varphi^{23\bar{3}},
\varphi^{12\overline{12}}, \varphi^{12\overline{13}}, \varphi^{12\overline{23}} \}=\C\{\varphi^{12\overline{12}}, \varphi^{12\overline{13}}, \varphi^{12\overline{23}}, \varphi^{13\overline{12}},
 \varphi^{23\overline{12}} \}\subset\ker\pb_{\phi(t)}.
\]
So
\[
\C\{\pb_{\phi(t)}\varphi^{13\bar{3}},\pb_{\phi(t)}\varphi^{23\bar{3}},
\varphi^{12\overline{12}}, \p\xi_5^{1,2}, \p\xi_6^{1,2} \}=\C\{\varphi^{12\overline{12}}, \varphi^{12\overline{13}}, \varphi^{12\overline{23}}, \varphi^{13\overline{12}},
 \varphi^{23\overline{12}} \}\subset\ker\pb_{\phi(t)}.
\]
Meanwhile, because the rank of
$\pb_{\phi(t)}:\C\{\varphi^{13\overline{13}}, \varphi^{13\overline{23}},\varphi^{23\overline{13}}, \varphi^{23\overline{23}}\}\longrightarrow\C\{\varphi^{12\overline{123}} \}$ is $1$,
we may choose $\xi_1^{2,2}, \xi_2^{2,2}, \xi_3^{2,2}, \xi_4^{2,2}\in\ker\pb_{\phi(t)}$ so that
\[
\C\{ \xi_1^{2,2}, \xi_2^{2,2}, \xi_3^{2,2}, \xi_4^{2,2}\}
=\C\{\varphi^{13\overline{13}}, \varphi^{13\overline{23}},\varphi^{23\overline{13}}, \varphi^{23\overline{23}}\}~\text{and}~
\pb_{\phi(t)}\xi_4^{2,2}=\varphi^{12\overline{123}}.
\]

For $\wedge^{1,3}$, we already see that $\pb_{\phi(t)}\xi_6^{1,2}\in \C\{\varphi^{1\overline{123}}, \varphi^{2\overline{123}}\}$ is nonzero, therefore we may choose $\xi_1^{1,3}$ so that
\[
\C\{\pb_{\phi(t)}\xi_6^{1,2}, \xi_1^{1,3}\}=\C\{\varphi^{1\overline{123}}, \varphi^{2\overline{123}}\}.
\]
\end{proof}
\newpage
By using Proposition \ref{prop-D(t)=0}, we immediately get the following diagrammatic view of the structure of $(\wedge^{\bullet}, d_{\phi(t)}=\p+\pb_{\phi(t)})$ for $(t_{11},t_{12},t_{21},t_{22})\neq 0, D(t)=0$:
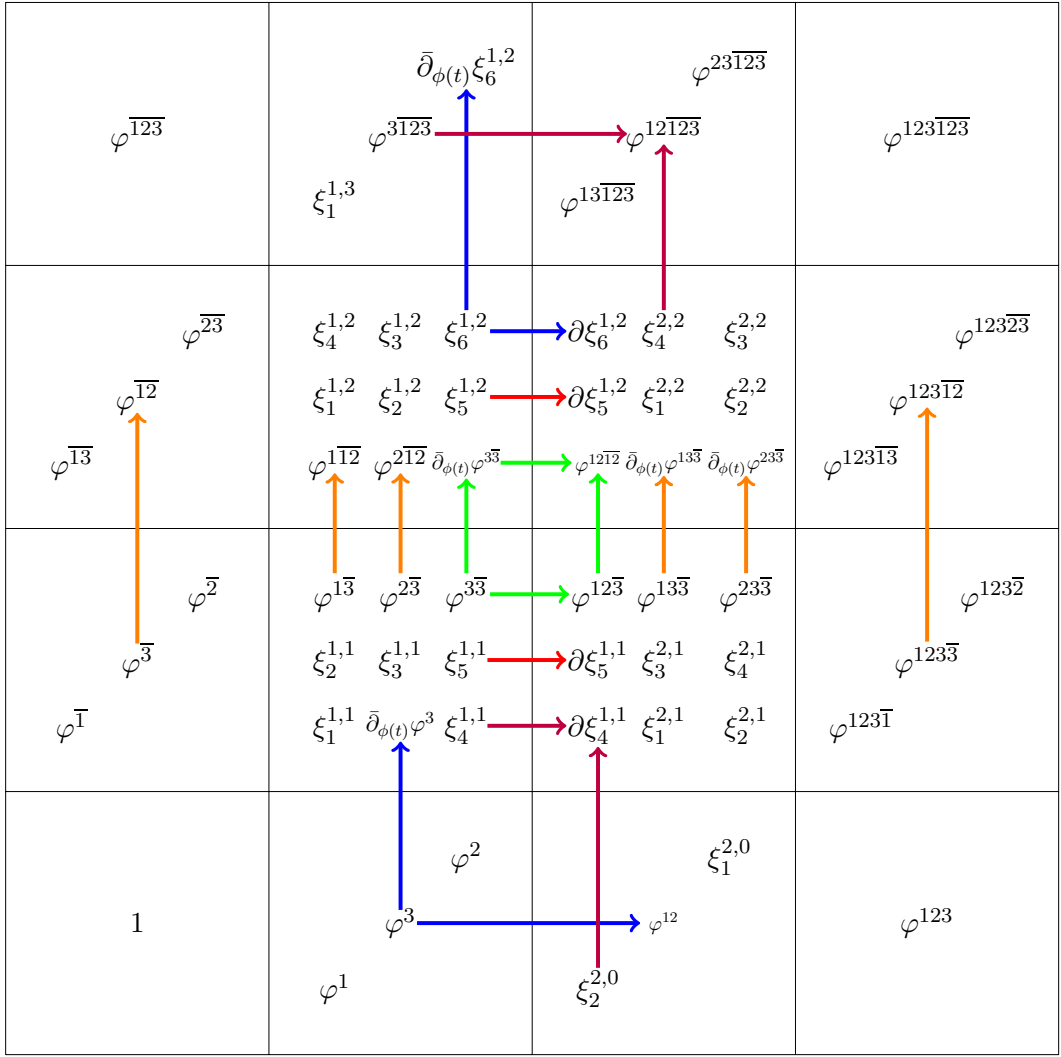
\begin{center}
\begin{tikzpicture}[scale=3.5]
\draw[help lines,black] (0,0) grid (4,4);
\draw [line width=1.5,->,blue] (3/2+1/16,1/2) -- (5/2-0.09,1/2);
\draw [line width=1.5,->,orange] (1/2,3/2+1/16) -- (1/2,5/2-0.06);
\draw [line width=1.5,->,blue] (3/2,1/2+0.05) -- (3/2,5/4-0.06);

\draw[line width=1.5,->,purple]  (9/4,1/4+0.08) -- (9/4,5/4-0.08);
\draw [line width=1.5,->,purple] (7/4+0.08,5/4) -- (9/4-0.12,5/4);
\draw [line width=1.5,->,red] (7/4+0.08,3/2) -- (9/4-0.12,3/2);
\draw [line width=1.5,->,orange] (5/4,7/4+0.08) -- (5/4,9/4-0.04);
\draw [line width=1.5,->,orange] (3/2,7/4+0.08) -- (3/2,9/4-0.04);
\draw [line width=1.5,->,green] (7/4+0.09,7/4) -- (9/4-0.11,7/4);
\draw [line width=1.5,->,green] (7/4,7/4+0.08) -- (7/4,9/4-0.06);
\draw [line width=1.5,->,green] (9/4,7/4+0.08) -- (9/4,9/4-0.04);
\draw [line width=1.5,->,green] (7/4+0.13,9/4) -- (9/4-0.1,9/4);

\draw [line width=1.5,->,red] (7/4+0.09,5/2) -- (9/4-0.12,5/2);
\draw [line width=1.5,->,blue] (7/4+0.09,11/4) -- (9/4-0.12,11/4);
\draw [line width=1.5,->,orange] (5/2,7/4+0.08) -- (5/2,9/4-0.05);
\draw [line width=1.5,->,orange] (45/16,7/4+0.08) -- (45/16,9/4-0.05);
\draw [line width=1.5,->,blue] (7/4,11/4+0.08) -- (7/4,15/4-0.08);

\draw [line width=1.5,->,orange] (7/2,3/2+0.07) -- (7/2,5/2-0.04);
\draw [line width=1.5,->,purple] (3/2+0.13,7/2) -- (5/2-0.14,7/2);
\draw [line width=1.5,->,purple] (5/2,11/4+0.08) -- (5/2,7/2-0.04);

\draw [fill] (1/2,1/2)   node {$1$};

\draw [fill=red] (5/4,1/4)  node {$\varphi^1$};
\draw [fill=red] (3/2,1/2) node {$\varphi^3$};
\draw [fill=red] (7/4,3/4) node {$\varphi^2$};

\draw [fill=green] (9/4,1/4) node {$\xi_2^{2,0}$};
\draw [fill=red] (5/2,1/2) node [scale=0.7]{$\varphi^{12}$};
\draw [fill=green] (11/4,3/4) node {$\xi_1^{2,0}$};

\draw [fill=blue] (7/2,1/2) node {$\varphi^{123}$};

\draw [fill=red] (1/4,5/4) node {$\varphi^{\overline{1}}$};
\draw [fill=red] (1/2,3/2) node {$\varphi^{\overline{3}}$};
\draw [fill=red] (3/4,7/4) node {$\varphi^{\overline{2}}$};

\draw [fill=green] (5/4,5/4) node {$\xi_1^{1,1}$};
\draw [fill=green] (3/2,5/4) node[scale=0.8] {$\pb_{\phi(t)}\varphi^3$};
\draw [fill=green] (7/4,5/4) node {$\xi_4^{1,1}$};
\draw [fill=green] (5/4,3/2) node {$\xi_2^{1,1}$};
\draw [fill=green] (3/2,3/2)  node {$\xi_3^{1,1}$};
\draw [fill=green] (7/4,3/2)  node {$\xi_5^{1,1}$};
\draw [fill=green] (5/4,7/4)  node {$\varphi^{1\overline{3}}$};
\draw [fill=green] (3/2,7/4)  node {$\varphi^{2\overline{3}}$};
\draw [fill=green] (7/4,7/4)  node {$\varphi^{3\overline{3}}$};

\draw [fill=green] (9/4,5/4) node {$\p\xi^{1,1}_4$};
\draw [fill=blue] (5/2,5/4) node {$\xi_1^{2,1}$};
\draw [fill=blue] (45/16,5/4) node {$\xi_2^{2,1}$};
\draw [fill=green] (9/4,3/2) node {$\p\xi_5^{1,1}$};
\draw [fill=blue] (5/2,3/2) node {$\xi_3^{2,1}$};
\draw [fill=blue] (45/16,3/2) node {$\xi_4^{2,1}$};
\draw [fill=green] (9/4,7/4) node {$\varphi^{12\overline{3}}$};
\draw [fill=blue] (5/2,7/4) node {$\varphi^{13\overline{3}}$};
\draw [fill=blue] (45/16,7/4) node {$\varphi^{23\overline{3}}$};

\draw [fill=orange] (13/4,5/4) node {$\varphi^{123\overline{1}}$};
\draw [fill=orange] (7/2,3/2) node {$\varphi^{123\overline{3}}$};
\draw [fill=orange] (15/4,7/4) node {$\varphi^{123\overline{2}}$};

\draw [fill=green] (1/4,9/4) node {$\varphi^{\overline{13}}$};
\draw [fill=red] (1/2,5/2) node {$\varphi^{\overline{12}}$};
\draw [fill=green] (3/4,11/4) node {$\varphi^{\overline{23}}$};

\draw [fill=green] (5/4,9/4) node {$\varphi^{1\overline{12}}$};
\draw [fill=green] (3/2,9/4) node {$\varphi^{2\overline{12}}$};
\draw [fill=green] (7/4,9/4) node[scale=0.7] {$\pb_{\phi(t)}\varphi^{3\overline{3}}$};
\draw [fill=blue] (5/4,5/2) node {$\xi_1^{1,2}$};
\draw [fill=blue] (3/2,5/2) node {$\xi_2^{1,2}$};
\draw [fill=blue] (7/4,5/2) node {$\xi_5^{1,2}$};
\draw [fill=blue] (5/4,11/4) node {$\xi_4^{1,2}$};
\draw [fill=blue] (3/2,11/4) node {$\xi_3^{1,2}$};
\draw [fill=blue] (7/4,11/4) node {$\xi_6^{1,2}$};

\draw [fill=green] (9/4,9/4) node[scale=0.7] {$\varphi^{12\overline{12}}$};
\draw [fill=blue] (5/2,9/4) node[scale=0.7] {$\pb_{\phi(t)}\varphi^{13\overline{3}}$};
\draw [fill=blue] (45/16,9/4) node[scale=0.7] {$\pb_{\phi(t)}\varphi^{23\overline{3}}$};
\draw [fill=blue] (9/4,5/2) node {$\p\xi_5^{1,2}$};
\draw [fill=orange] (5/2,5/2) node {$\xi_1^{2,2}$};
\draw [fill=orange] (45/16,5/2) node {$\xi_2^{2,2}$};
\draw [fill=blue] (9/4,11/4) node {$\p\xi_6^{1,2}$};
\draw [fill=orange] (5/2,11/4) node {$\xi_4^{2,2}$};
\draw [fill=orange] (45/16,11/4) node {$\xi_3^{2,2}$};

\draw [fill=magenta] (13/4,9/4) node {$\varphi^{123\overline{13}}$};
\draw [fill=orange] (7/2,5/2) node {$\varphi^{123\overline{12}}$};
\draw [fill=magenta] (15/4,11/4) node {$\varphi^{123\overline{23}}$};

\draw [fill=blue] (1/2,7/2) node {$\varphi^{\overline{123}}$};

\draw [fill=orange] (5/4,13/4) node {$\xi_1^{1,3}$};
\draw [fill=orange] (3/2,7/2) node {$\varphi^{3\overline{123}}$};
\draw [fill=orange] (7/4,15/4) node {$\pb_{\phi(t)}\xi_6^{1,2}$};

\draw [fill=magenta] (9/4,13/4) node {$\varphi^{13\overline{123}}$};
\draw [fill=orange] (5/2,7/2) node {$\varphi^{12\overline{123}}$};
\draw [fill=magenta] (11/4,15/4) node {$\varphi^{23\overline{123}}$};

\draw [fill] (7/2,7/2) node {$\varphi^{123\overline{123}}$};
\end{tikzpicture}
\end{center}

\subsection{The case that $(t_{11},t_{12},t_{21},t_{22})\neq 0, D(t)\neq0$}\label{subsec-(iii)}
\begin{proposition}\label{prop-D(t)not0}
There exist the following basis
\begin{itemize}
  \item $\wedge^{1,0}=\C\{ \underbrace{\varphi^{1}, \varphi^{2}}_{\ker\pb_{\phi(t)}}, \varphi^{3} \}$;
  \item $\wedge^{2,0}=\C\{\underbrace{\varphi^{12}}_{\ker\pb_{\phi(t)}}, \varphi^{13}, \varphi^{23}\}$;
  \item $\wedge^{1,1}=\C\{\underbrace{\pb_{\phi(t)}\varphi^{3}, \xi_1^{1,1}, \xi_2^{1,1}, \xi_3^{1,1}, \eta_1^{1,1}, \eta_2^{1,1}}_{\ker\pb_{\phi(t)}}, \varphi^{1\bar{3}}, \varphi^{2\bar{3}}, \varphi^{3\bar{3}} \}$;
  \item $\wedge^{2,1}=\C\{\underbrace{ \p\eta_1^{1,1}=\pb_{\phi(t)}\varphi^{23}, \p\eta_2^{1,1}=\pb_{\phi(t)}\varphi^{13}, \xi_2^{2,1}, \xi_3^{2,1}, \xi_4^{2,1}}_{\ker\pb_{\phi(t)}}, \varphi^{12\bar{3}}, \varphi^{13\bar{3}}, \varphi^{23\bar{3}} \}$;
  \item $\wedge^{1,2}=\C\{\underbrace{\pb_{\phi(t)}\varphi^{3\bar{3}}, \varphi^{1\overline{12}}, \varphi^{2\overline{12}}, \xi_1^{1,2}, \xi_2^{1,2}, \xi_3^{1,2}, \xi_4^{1,2}}_{\ker\pb_{\phi(t)}}, \varphi^{3\overline{13}}, \varphi^{3\overline{23}} \}$;
  \item $\wedge^{2,2}=\C\{\underbrace{\pb_{\phi(t)}\varphi^{13\bar{3}}, \pb_{\phi(t)}\varphi^{23\bar{3}}, \varphi^{12\overline{12}}, \varphi^{12\overline{13}}, \varphi^{12\overline{23}}, \xi_1^{2,2}, \xi_2^{2,2}, \xi_3^{2,2}}_{\ker\pb_{\phi(t)}}, \xi_4^{2,2} \}$;
  \item $\wedge^{1,3}=\C\{\underbrace{\pb_{\phi(t)}\varphi^{3\overline{13}}, \pb_{\phi(t)}\varphi^{3\overline{23}}}_{\in\ker\pb_{\phi(t)}}, \varphi^{3\overline{123}} \}$,
\end{itemize}
such that
\begin{itemize}
  \item $\mathbb{C}\{\xi_1^{1,1},\xi_2^{1,1},\xi_3^{1,1},\pb_{\phi(t)}\varphi^3\}
      =\mathbb{C}\{\varphi^{1\bar{1}},\varphi^{1\bar{2}},
      \varphi^{2\bar{1}},\varphi^{2\bar{2}}\}$ and
  \begin{align*}
   \mathbb{C}\{\eta_1^{1,1},\eta_2^{1,1}\}&=
   \mathbb{C}\{\varphi^{3\bar{1}},\varphi^{3\bar{2}}\},\quad \p\eta_1^{1,1}=\pb_{\phi(t)}\varphi^{23},~\p\eta_2^{1,1}=\pb_{\phi(t)}\varphi^{13},\\
   \p\varphi^{3\bar{3}}&=-\varphi^{12\bar{3}},~\pb_{\phi(t)}\varphi^{3\bar{3}}=\varphi^{1\overline{12}},~\pb_{\phi(t)}\varphi^{2\bar{3}}=\varphi^{2\overline{12}};
  \end{align*}
  \item $\mathbb{C}\{\varphi^{12\bar{3}},\xi_1^{2,1},\xi_2^{2,1},\xi_3^{2,1},
      \xi_4^{2,1}\}=\mathbb{C}\{\varphi^{12\bar{3}},\varphi^{13\bar{1}},
      \varphi^{13\bar{2}},\varphi^{23\bar{1}},\varphi^{23\bar{2}}\}$ and
      \[
      \pb_{\phi(t)}\varphi^{12\bar{3}}=-\varphi^{12\overline{12}},
      ~ \C\{\pb_{\phi(t)}\varphi^{13\overline{3}},\pb_{\phi(t)}\varphi^{23\overline{3}}\}
      \subset\C\{\varphi^{13\overline{12}},\varphi^{23\overline{12}},
      \varphi^{12\overline{13}},\varphi^{12\overline{23}}\}~\text{is~of~dimension}~2;
      \]
  \item $\mathbb{C}\{\xi_1^{1,2},\xi_2^{1,2},\xi_3^{1,2},\xi_4^{1,2},
      \pb_{\phi(t)}\varphi^{3\overline{3}}\}=
      \mathbb{C}\{\varphi^{1\overline{13}},\varphi^{1\overline{23}},
      \varphi^{2\overline{13}},\varphi^{2\overline{23}},\varphi^{3\overline{12}}\}$
      and
      \[
\mathbb{C}\{\pb_{\phi(t)}\varphi^{3\overline{13}}, \pb_{\phi(t)}\varphi^{3\overline{23}}\}=\C\{\varphi^{1\overline{123}}, \varphi^{2\overline{123}}\},\quad \text{with}~\p\varphi^{3\overline{13}}=-\varphi^{12\overline{13}}, ~\p\varphi^{3\overline{23}}=-\varphi^{12\overline{23}};
      \]
  \item $\mathbb{C}\{\xi_1^{2,2},\xi_2^{2,2},\xi_3^{2,2},\xi_4^{2,2}\}
      =\mathbb{C}\{\varphi^{13\overline{13}},\varphi^{13\overline{23}},
      \varphi^{23\overline{13}},\varphi^{23\overline{23}}\}$ with
      \[
      \xi_1^{2,2},\xi_2^{2,2},\xi_3^{2,2}\in\ker\pb_{\phi(t)}~\text{and}
      ~\pb_{\phi(t)}\xi_4^{2,2}=\varphi^{12\overline{123}};
      \]
  \item $\C\{\pb_{\phi(t)}\varphi^{3\overline{13}}, \pb_{\phi(t)}\varphi^{3\overline{23}} \}=
      \mathbb{C}\{\varphi^{1\overline{123}},\varphi^{2\overline{123}}\}$.
\end{itemize}
\end{proposition}

\begin{proof}From Subsection \ref{subsubsec-10to11}, we know that $\ker\pb_{\phi(t)}\cap\wedge^{1,0}=\C\{\varphi^{1}, \varphi^{2}\}$ and the rank of $\pb_{\phi(t)}:\C\{\varphi^3\}\to
\C\{\varphi^{1\bar{1}},\varphi^{2\bar{1}},\varphi^{1\bar{2}},\varphi^{2\bar{2}}\}$ is $1$.

From Subsection \ref{subsubsec-20to21}, we know that $\dim\ker\pb_{\phi(t)}\cap\wedge^{2,0}=1$ and $\pb_{\phi(t)}\varphi^{12}=0$. On the other hand, the rank of $\pb_{\phi(t)}:\C\{\varphi^{13}, \varphi^{23}\}\to
\C\{\varphi^{12\bar{1}},\varphi^{12\bar{2}}\}$ is $2$. So we have
\[
\C\{\pb_{\phi(t)}\varphi^{12\bar{1}}, \pb_{\phi(t)}\varphi^{12\bar{2}}\}
=
\C\{\varphi^{12\bar{1}},\varphi^{12\bar{2}}\}.
\]

From Subsection \ref{subsubsec-11to12}, we know that $\dim\ker\pb_{\phi(t)}\cap\wedge^{1,1}=6$ and the rank of
\[
\pb_{\phi(t)}:\C\{\varphi^{1\bar{3}}, \varphi^{2\bar{3}}, \varphi^{3\bar{1}}, \varphi^{3\bar{2}}\}\longrightarrow\C\{\varphi^{1\overline{12}}, \varphi^{2\overline{12}}\},
\]
is equal to $\text{rank}(T)=2$. So we may choose $\xi_1^{1,1},\xi_2^{1,1},\xi_3^{1,1}$ so that $\mathbb{C}\{\xi_1^{1,1},\xi_2^{1,1},\xi_3^{1,1},\pb_{\phi(t)}\varphi^3\}
      =\mathbb{C}\{\varphi^{1\bar{1}},\varphi^{1\bar{2}},
      \varphi^{2\bar{1}},\varphi^{2\bar{2}}\}$,
and $\eta_1^{1,1},\eta_2^{1,1}$ so that (indeed, we may set $\eta_1^{1,1}=-t_{11}\varphi^{3\bar{1}}-t_{12}\varphi^{3\bar{1}}$ and $\eta_2^{1,1}=-t_{21}\varphi^{3\bar{1}}-t_{22}\varphi^{3\bar{1}}$)
\[
\mathbb{C}\{\eta_1^{1,1},\eta_2^{1,1}\}=
   \mathbb{C}\{\varphi^{3\bar{1}},\varphi^{3\bar{2}}\},\quad \p\eta_1^{1,1}=\pb_{\phi(t)}\varphi^{23},~\p\eta_2^{1,1}=\pb_{\phi(t)}\varphi^{13}.
\]
\[
\ker\pb_{\phi(t)}=\C\{\varphi^{1\bar{1}},\varphi^{1\bar{2}},
\varphi^{2\bar{1}},\varphi^{2\bar{2}},\xi_4^{1,1},\xi_5^{1,1}\}.
\]
On the other hand, $\pb_{\phi(t)}\varphi^{1\bar{3}}=\varphi^{1\overline{12}}, \pb_{\phi(t)}\varphi^{2\bar{3}}=\varphi^{2\overline{12}}$ and the rank of
$\pb_{\phi(t)}:\C\{\varphi^{3\bar{3}}\}\longrightarrow\C\{\varphi^{1\overline{13}}, \varphi^{1\overline{23}}, \varphi^{2\overline{13}}, \varphi^{2\overline{23}}, \varphi^{3\overline{12}}\}$ is $1$.
Notice that
\[
\p\pb_{\phi(t)}\varphi^{3\bar{3}}=-\varphi^{12\overline{12}}
=-\pb_{\phi(t)}\varphi^{12\bar{3}}=\pb_{\phi(t)}\p\varphi^{3\bar{3}},
\]

From Subsection \ref{subsubsec-21to22}, we know that $\mathbb{C}\{\p\xi_4^{1,1},\p\xi_5^{1,1}\}=\C\{\varphi^{12\bar{1}},\varphi^{12\bar{2}}\}
\subset\ker\pb_{\phi(t)}$, the rank of $\pb_{\phi(t)}:\C\{\varphi^{12\bar{3}}, \varphi^{13\bar{1}}, \varphi^{13\bar{2}}, \varphi^{23\bar{1}},\varphi^{23\bar{2}}\}\longrightarrow\C\{\varphi^{12\overline{12}}\}$ is $1$, $\pb_{\phi(t)}\varphi^{13\bar{3}}$ and $\pb_{\phi(t)}\varphi^{23\bar{3}}$ are linearly independent. We may choose $\xi_1^{2,1},\xi_2^{2,1},\xi_3^{2,1},\xi_4^{2,1}\in\ker\pb_{\phi(t)}$ so that
\[
\mathbb{C}\{\varphi^{12\bar{3}},\xi_1^{2,1},\xi_2^{2,1},\xi_3^{2,1},\xi_4^{2,1}\}
=\mathbb{C}\{\varphi^{12\bar{3}},\varphi^{13\bar{1}},
\varphi^{13\bar{2}},\varphi^{23\bar{1}},\varphi^{23\bar{2}}\}.
\]

For $\wedge^{1,2}$, we already see that $\varphi^{1\overline{12}}=\pb_{\phi(t)}\varphi^{1\bar{3}}, \varphi^{2\overline{12}}=\pb_{\phi(t)}\varphi^{2\bar{3}}$ and the rank of
$
\pb_{\phi(t)}:\C\{\varphi^{3\bar{3}}\}\longrightarrow\C\{\varphi^{1\overline{13}}, \varphi^{1\overline{23}}, \varphi^{2\overline{13}}, \varphi^{2\overline{23}}, \varphi^{3\overline{12}}\}
$ is $1$. From Subsection \ref{subsubsec-12to13}, we know that $\varphi^{1\overline{12}}, \varphi^{1\overline{13}}, \varphi^{1\overline{23}}, \varphi^{2\overline{12}}, \varphi^{2\overline{13}}, \varphi^{2\overline{23}}, \varphi^{3\overline{12}}$ are in $\ker\pb_{\phi(t)}$ and the rank of $\pb_{\phi(t)}:\C\{\varphi^{3\overline{13}}, \varphi^{3\overline{23}}\}\longrightarrow\C\{\varphi^{1\overline{123}}, \varphi^{2\overline{123}}\}$ is $1$. So we may choose $\xi_1^{1,2},\cdots,\xi_6^{1,2}$ so that
$\mathbb{C}\{\xi_1^{1,2},\xi_2^{1,2},\xi_3^{1,2},\xi_4^{1,2},
\pb_{\phi(t)}\varphi^{3\overline{3}}\}=
\mathbb{C}\{\varphi^{1\overline{13}},\varphi^{1\overline{23}},
\varphi^{2\overline{13}},\varphi^{2\overline{23}},\varphi^{3\overline{12}}\}$
and
\[
\mathbb{C}\{\xi_5^{1,2}, \xi_6^{1,2}\}=\C\{\varphi^{3\overline{13}}, \varphi^{3\overline{23}}\},~\text{with}~
\pb_{\phi(t)}\xi_5^{1,2}=0,~\pb_{\phi(t)}\xi_6^{1,2}\neq0.
\]

For $\wedge^{2,2}$, we already see that
\[
\pb_{\phi(t)}\varphi^{13\bar{3}},\pb_{\phi(t)}\varphi^{23\bar{3}}\in\C
\{\varphi^{13\overline{12}},\varphi^{12\overline{13}},\varphi^{12\overline{23}},
\varphi^{23\overline{12}}\},
\]
are linearly independent. From Subsection \ref{subsubsec-22to23}, we know that 
\[
\varphi^{12\overline{12}}, \varphi^{12\overline{13}}, \varphi^{12\overline{23}}, \varphi^{13\overline{12}}, \varphi^{23\overline{12}}\in\ker\pb_{\phi(t)}
\]
and the rank of
$\pb_{\phi(t)}:\C\{\varphi^{13\overline{13}}, \varphi^{13\overline{23}},\varphi^{23\overline{13}}, \varphi^{23\overline{23}}\}\longrightarrow\C\{\varphi^{12\overline{123}}\}$ is $1$. Notice that $\C\{\varphi^{12\overline{13}},\varphi^{12\overline{23}} \}=
\mathbb{C}\{\p\xi_5^{1,2}, \p\xi_6^{1,2}\}$ and
\[
\C\{\pb_{\phi(t)}\varphi^{13\bar{3}},\pb_{\phi(t)}\varphi^{23\bar{3}},
\varphi^{12\overline{12}}, \varphi^{12\overline{13}}, \varphi^{12\overline{23}} \}=\C\{\varphi^{12\overline{12}}, \varphi^{12\overline{13}}, \varphi^{12\overline{23}}, \varphi^{13\overline{12}},
 \varphi^{23\overline{12}} \}\subset\ker\pb_{\phi(t)}.
\]
So
\[
\C\{\pb_{\phi(t)}\varphi^{13\bar{3}},\pb_{\phi(t)}\varphi^{23\bar{3}},
\varphi^{12\overline{12}}, \p\xi_5^{1,2}, \p\xi_6^{1,2} \}=\C\{\varphi^{12\overline{12}}, \varphi^{12\overline{13}}, \varphi^{12\overline{23}}, \varphi^{13\overline{12}},
 \varphi^{23\overline{12}} \}\subset\ker\pb_{\phi(t)}.
\]
Meanwhile, because the rank of
$\pb_{\phi(t)}:\C\{\varphi^{13\overline{13}}, \varphi^{13\overline{23}},\varphi^{23\overline{13}}, \varphi^{23\overline{23}}\}\longrightarrow\C\{\varphi^{12\overline{123}} \}$ is $1$,
we may choose $\xi_1^{2,2}, \xi_2^{2,2}, \xi_3^{2,2}, \xi_4^{2,2}\in\ker\pb_{\phi(t)}$ so that
\[
\C\{ \xi_1^{2,2}, \xi_2^{2,2}, \xi_3^{2,2}, \xi_4^{2,2}\}
=\C\{\varphi^{13\overline{13}}, \varphi^{13\overline{23}},\varphi^{23\overline{13}}, \varphi^{23\overline{23}}\}~\text{and}~
\pb_{\phi(t)}\xi_4^{2,2}=\varphi^{12\overline{123}}.
\]

For $\wedge^{1,3}$, we already see that $\C\{\pb_{\phi(t)}\varphi^{3\overline{13}}, \pb_{\phi(t)}\varphi^{3\overline{23}} \}=
\mathbb{C}\{\varphi^{1\overline{123}},\varphi^{2\overline{123}}\}$ and it is obvious that $\pb_{\phi(t)}\varphi^{3\overline{123}}=0$.
\end{proof}
\newpage
By using Proposition \ref{prop-D(t)not0}, we immediately get the following diagrammatic view of the structure of $(\wedge^{\bullet}, d_{\phi(t)}=\p+\pb_{\phi(t)})$ for $(t_{11},t_{12},t_{21},t_{22})\neq 0, D(t)\neq0$:

\begin{figure}[ht]
\begin{center}
\begin{tikzpicture}[scale=3.5]

\draw[help lines,black] (0,0) grid (4,4);
\draw [line width=1.5,->,blue] (3/2+1/16,1/2) -- (5/2-0.09,1/2);
\draw [line width=1.5,->,orange] (1/2,3/2+1/16) -- (1/2,5/2-0.06);
\draw [line width=1.5,->,blue] (3/2,1/2+0.05) -- (3/2,5/4-0.07);

\draw[line width=1.5,->,purple]  (9/4,1/4+0.06) -- (9/4,5/4-0.06);
\draw [line width=1.5,->,purple] (7/4+0.06,5/4) -- (9/4-0.11,5/4);
\draw [line width=1.5,->,purple] (7/4+0.06,3/2) -- (11/4-0.11,3/2);
\draw [line width=1.5,->,orange] (5/4,7/4+0.06) -- (5/4,9/4-0.04);
\draw [line width=1.5,->,orange] (3/2,7/4+0.06) -- (3/2,9/4-0.04);
\draw [line width=1.5,->,green] (7/4+0.08,7/4) -- (9/4-0.1,7/4);
\draw [line width=1.5,->,green] (7/4,7/4+0.06) -- (7/4,9/4-0.05);
\draw [line width=1.5,->,green] (9/4,7/4+0.06) -- (9/4,9/4-0.03);
\draw [line width=1.5,->,green] (7/4+0.13,9/4) -- (9/4-0.09,9/4);
\draw [line width=1.5,->,purple] (11/4,3/4+0.06) -- (11/4,3/2-0.06);

\draw [line width=1.5,->,blue] (5/4+0.11,5/2) -- (9/4-0.12,5/2);
\draw [line width=1.5,->,blue] (7/4+0.11,11/4) -- (9/4-0.12,11/4);
\draw [line width=1.5,->,orange] (5/2,7/4+0.06) -- (5/2,9/4-0.05);
\draw [line width=1.5,->,orange] (45/16,7/4+0.06) -- (45/16,9/4-0.05);
\draw [line width=1.5,->,blue] (7/4,11/4+0.07) -- (7/4,15/4-0.07);
\draw [line width=1.5,->,blue] (5/4,5/2+0.07) -- (5/4,13/4-0.07);

\draw [line width=1.5,->,orange] (7/2,3/2+0.06) -- (7/2,5/2-0.03);
\draw [line width=1.5,->,purple] (3/2+0.13,7/2) -- (5/2-0.15,7/2);
\draw [line width=1.5,->,purple] (5/2,11/4+0.08) -- (5/2,7/2-0.03);

\draw [fill] (1/2,1/2)   node {$1$};

\draw [fill=red] (5/4,1/4)  node {$\varphi^1$};
\draw [fill=red] (3/2,1/2) node {$\varphi^3$};
\draw [fill=red] (7/4,3/4) node {$\varphi^2$};

\draw [fill=green] (9/4,1/4) node {$\varphi^{23}$};
\draw [fill=red] (5/2,1/2) node {$\varphi^{12}$};
\draw [fill=green] (11/4,3/4) node {$\varphi^{13}$};

\draw [fill=blue] (7/2,1/2) node {$\varphi^{123}$};

\draw [fill=red] (1/4,5/4) node {$\varphi^{\overline{1}}$};
\draw [fill=red] (1/2,3/2) node {$\varphi^{\overline{3}}$};
\draw [fill=red] (3/4,7/4) node {$\varphi^{\overline{2}}$};

\draw [fill=green] (5/4,5/4) node {$\xi_1^{1,1}$};
\draw [fill=green] (3/2,5/4) node [scale=0.8]{$\pb_{\phi(t)}\varphi^3$};
\draw [fill=green] (7/4,5/4) node {$\eta_1^{1,1}$};
\draw [fill=green] (5/4,3/2) node {$\xi_2^{1,1}$};
\draw [fill=green] (3/2,3/2)  node {$\xi_3^{1,1}$};
\draw [fill=green] (7/4,3/2)  node {$\eta_2^{1,1}$};
\draw [fill=green] (5/4,7/4)  node {$\varphi^{1\overline{3}}$};
\draw [fill=green] (3/2,7/4)  node {$\varphi^{2\overline{3}}$};
\draw [fill=green] (7/4,7/4)  node {$\varphi^{3\overline{3}}$};

\draw [fill=green] (9/4+0.04,5/4) node [scale=0.8]{$\pb_{\phi(t)}\varphi^{23}$};
\draw [fill=blue] (5/2+0.02,5/4) node {$\xi_1^{2,1}$};
\draw [fill=blue] (45/16,5/4) node {$\xi_2^{2,1}$};
\draw [fill=green] (9/4,3/2) node {$\xi_4^{2,1}$};
\draw [fill=blue] (5/2+0.02,3/2) node {$\xi_3^{2,1}$};
\draw [fill=blue] (45/16,3/2) node [scale=0.8]{$\pb_{\phi(t)}\varphi^{13}$};
\draw [fill=green] (9/4,7/4) node {$\varphi^{12\overline{3}}$};
\draw [fill=blue] (5/2+0.02,7/4) node {$\varphi^{13\overline{3}}$};
\draw [fill=blue] (45/16,7/4) node {$\varphi^{23\overline{3}}$};

\draw [fill=orange] (13/4,5/4) node {$\varphi^{123\overline{1}}$};
\draw [fill=orange] (7/2,3/2) node {$\varphi^{123\overline{3}}$};
\draw [fill=orange] (15/4,7/4) node {$\varphi^{23\overline{2}}$};

\draw [fill=green] (1/4,9/4) node {$\varphi^{\overline{13}}$};
\draw [fill=red] (1/2,5/2) node {$\varphi^{\overline{12}}$};
\draw [fill=green] (3/4,11/4) node {$\varphi^{\overline{23}}$};

\draw [fill=green] (5/4,9/4) node {$\varphi^{1\overline{12}}$};
\draw [fill=green] (3/2,9/4) node {$\varphi^{2\overline{12}}$};
\draw [fill=green] (7/4,9/4) node [scale=0.7]{$\pb_{\phi(t)}\varphi^{3\overline{3}}$};
\draw [fill=blue] (5/4,5/2) node {$\varphi^{3\overline{13}}$};
\draw [fill=blue] (3/2,5/2) node {$\xi_2^{1,2}$};
\draw [fill=blue] (7/4,5/2) node {$\xi_1^{1,2}$};
\draw [fill=blue] (5/4,11/4) node {$\xi_4^{1,2}$};
\draw [fill=blue] (3/2,11/4) node {$\xi_3^{1,2}$};
\draw [fill=blue] (7/4-0.04,11/4) node {$\varphi^{3\overline{23}}$};

\draw [fill=green] (9/4,9/4) node [scale=0.7]{$\varphi^{12\overline{12}}$};
\draw [fill=blue] (5/2,9/4) node [scale=0.7]{$\pb_{\phi(t)}\varphi^{13\overline{3}}$};
\draw [fill=blue] (45/16,9/4) node [scale=0.7]{$\pb_{\phi(t)}\varphi^{23\overline{3}}$};
\draw [fill=blue] (9/4,5/2) node {$\varphi^{12\overline{13}}$};
\draw [fill=orange] (5/2,5/2) node {$\xi_1^{2,2}$};
\draw [fill=orange] (45/16,5/2) node {$\xi_2^{2,2}$};
\draw [fill=blue] (9/4,11/4) node {$\varphi^{12\overline{23}}$};
\draw [fill=orange] (5/2,11/4) node {$\xi_4^{2,2}$};
\draw [fill=orange] (45/16,11/4) node {$\xi_3^{2,2}$};

\draw [fill=magenta] (13/4,9/4) node {$\varphi^{123\overline{13}}$};
\draw [fill=orange] (7/2,5/2) node {$\varphi^{123\overline{12}}$};
\draw [fill=magenta] (15/4,11/4) node {$\varphi^{12\overline{23}}$};

\draw [fill=blue] (1/2,7/2) node {$\varphi^{\overline{123}}$};

\draw [fill=orange] (5/4,13/4) node {$\pb_{\phi(t)}\varphi^{3\overline{13}}$};
\draw [fill=orange] (3/2,7/2) node {$\varphi^{3\overline{123}}$};
\draw [fill=orange] (7/4,15/4) node {$\pb_{\phi(t)}\varphi^{3\overline{23}}$};

\draw [fill=magenta] (9/4,13/4) node {$\varphi^{13\overline{123}}$};
\draw [fill=orange] (5/2,7/2) node {$\varphi^{12\overline{123}}$};
\draw [fill=magenta] (11/4,15/4) node {$\varphi^{23\overline{123}}$};

\draw [fill] (7/2,7/2) node {$\varphi^{123\overline{123}}$};
\end{tikzpicture}
\end{center}
\end{figure}

\subsection{Proof of Theorem \ref{thm-main-0}}
\begin{proof}[Proof of Theorem \ref{thm-main-0}]
The conclusions of Theorem \ref{thm-main-0} follows immediately from the two diagrams in Subsection \ref{subsec-(ii)} and \ref{subsec-(iii)}.
\end{proof}
\section{The dimension of cohomologies and $E_r^{p,q}(X_t)$}\label{sec-5}
In this section, we compute the Dolbeault/Bott-Chern cohomology and the Fr\"olicher spectral sequence $(E_r^{\bullet,\bullet}, d_r^{\bullet,\bullet})$ of the double complex $(\wedge^{\bullet}, d_{\phi(t)}=\p+\pb_{\phi(t)})$ by counting zigzags in diagrams obtained in Subsection \ref{subsec-(ii)} and \ref{subsec-(iii)}.

According to Theorem \ref{thm-cohomology-Er} and the fact that $(\wedge^{\bullet}, d_{\phi(t)}=\p+\pb_{\phi(t)})\hookrightarrow (A^{\bullet}(X), d_{\phi(t)}=\p+\pb_{\phi(t)})$ is an $E_1$-isomorphism, we immediately get the following

\renewcommand\arraystretch{1.5}
\begin{table}[!htbp]
\caption{Dimensions of $H^{\bullet,\bullet}_{\pb_{\phi(t)}}(X), H^{\bullet,\bullet}_{BC\phi(t)}(X)$ and $\im d_1^{\bullet,\bullet}(X_t)$}
\centering
\begin{center}
\begin{tabular}{|c|c|c|c|c|c|c|c|c|c|c|c|}
\hline
& $(i)$ & $(ii)$ & $(iii)$ & & $(i)$ & $(ii)$ & $(iii)$& & $(i)$ & $(ii)$ & $(iii)$ \\
\hline
 $h^{1,0}_{\pb_{\phi(t)}}$ & $3$ & $2$ & $2$ &$h^{1,0}_{BC,\phi(t)}$& $2$ & $2$ & $2$ & $\im d_1^{1,0}$& $1$ & $0$ & $0$ \\
\hline
 $h^{0,1}_{\pb_{\phi(t)}}$ & $2$ & $2$ & $2$ &$h^{0,1}_{BC,\phi(t)}$& $2$ & $2$ & $2$ &$\im d_1^{0,1}$&  $0$ & $0$ & $0$ \\
\hline
 $h^{2,0}_{\pb_{\phi(t)}}$ & $3$ & $2$ & $1$ & $h^{2,0}_{BC,\phi(t)}$& $3$ & $2$ & $1$ & $\im d_1^{2,0}$&$0$ & $0$ & $0$ \\
\hline
 $h^{1,1}_{\pb_{\phi(t)}}$ & $6$ & $5$ & $5$ & $h^{1,1}_{BC,\phi(t)}$& $4$ & $4$ & $4$ & $\im d_1^{1,1}$& $2$ & $1$ & $0$ \\
\hline
 $h^{0,2}_{\pb_{\phi(t)}}$ & $2$ & $2$ & $2$ &$h^{0,2}_{BC,\phi(t)}$ & $3$ & $3$ & $3$ &$\im d_1^{0,2}$ & $0$ & $0$ & $0$ \\
\hline
 $h^{3,0}_{\pb_{\phi(t)}}$ & $1$ & $1$ & $1$ &$h^{3,0}_{BC,\phi(t)}$ & $1$ & $1$ & $1$ &$\im d_1^{3,0}$ & $0$ & $0$ & $0$  \\
\hline
$h^{2,1}_{\pb_{\phi(t)}}$ & $6$ & $5$ & $4$ &$h^{2,1}_{BC,\phi(t)}$ & $6$ & $6$ & $6$ &$\im d_1^{2,1}$ & $0$ & $0$ & $0$  \\
\hline
$h^{1,2}_{\pb_{\phi(t)}}$ & $6$ & $5$ & $4$ &$h^{1,2}_{BC,\phi(t)}$ & $6$ & $6$ & $6$ &$\im d_1^{1,2}$ & $2$ & $1$ & $0$  \\
\hline
$h^{0,3}_{\pb_{\phi(t)}}$ & $1$ & $1$ & $1$ &$h^{0,3}_{BC,\phi(t)}$ & $1$ & $1$ & $1$ &$\im d_1^{0,3}$ & $0$ & $0$ & $0$  \\
\hline
$h^{3,1}_{\pb_{\phi(t)}}$ & $2$ & $2$ & $2$ &$h^{3,1}_{BC,\phi(t)}$ & $2$ & $2$ & $2$ &$\im d_1^{3,1}$ & $0$ & $0$ & $0$  \\
\hline
$h^{2,2}_{\pb_{\phi(t)}}$ & $6$ & $5$ & $5$ &$h^{2,2}_{BC,\phi(t)}$ & $8$ & $7$ & $7$ &$\im d_1^{2,2}$ & $0$ & $0$ & $0$  \\
\hline
$h^{1,3}_{\pb_{\phi(t)}}$ & $3$ & $2$ & $1$ &$h^{1,3}_{BC,\phi(t)}$ & $2$ & $2$ & $2$ &$\im d_1^{1,3}$ & $1$ & $0$ & $0$  \\
\hline
$h^{3,2}_{\pb_{\phi(t)}}$ & $2$ & $2$ & $2$ &$h^{3,2}_{BC,\phi(t)}$ & $3$ & $3$ & $3$ &$\im d_1^{3,2}$ & $0$ & $0$ & $0$  \\
\hline
$h^{2,3}_{\pb_{\phi(t)}}$ & $3$ & $2$ & $2$ &$h^{2,3}_{BC,\phi(t)}$ & $3$ & $3$ & $3$ &$\im d_1^{2,3}$ & $0$ & $0$ & $0$  \\
\hline
\end{tabular}
\end{center}
\end{table}
\noindent where $\im d_1^{p,q}$ is short for $\dim\im d_1^{p,q}(X_t)=\dim\im d_1^{p,q}(X,d_{\phi(t)})$ and we know that $H^{\bullet,\bullet}_{\pb_{\phi(t)}}(X)$ is isomorphic to $H^{\bullet,\bullet}_{\pb_t}(X_t)$ \cite{Xia19dDol}. It is clear that for $r>1$, $d_r^{p,q}=0$ for any $(p,q)$ because all zigzags appeared have length less than $4$. Furthermore, because
\[
E_{r+1}^{p,q}\cong \frac{\ker d_r^{p,q}}{\im d_r^{p-r,q+r-1}}\quad\text{and}\quad E_r^{p,q}/\ker d_r^{p,q}\cong\im d_r^{p,q},\quad\forall r\geq 1,
\]
we have
\begin{equation}\label{eq-dimE2tor}
      \left\{
      \begin{array}{ll}
\dim E_2^{p,q}&=\dim E_1^{p,q}-\dim \im d_1^{p,q}-\dim \im d_1^{p-1,q},\\
\dim E_3^{p,q}&=\dim E_2^{p,q}-\dim \im d_2^{p,q}-\dim \im d_2^{p-2,q+1},\\
&\vdots\\
\dim E_{r+1}^{p,q}&=\dim E_r^{p,q}-\dim \im d_r^{p,q}-\dim \im d_r^{p-r,q+r-1}.\\
\end{array} \right.
\end{equation}
It follows that
\renewcommand\arraystretch{1.5}
\begin{table}[!htbp]
\caption{$e_r^{\bullet,\bullet}=\dim E_r^{\bullet,\bullet}(X_t)$ and $\chi^p_q(E_2)=\sum_{i=0}^q (-1)^{q-i}e_2^{p,i}$}
\centering
\begin{center}
\begin{tabular}{|c|c|c|c|c|c|c|c|c|c|c|c|}
\hline
 & $(i)$ & $(ii)$ & $(iii)$ & & $(i)$ & $(ii)$ & $(iii)$ & & $(i)$ & $(ii)$ & $(iii)$ \\
\hline
$e_1^{1,0}$ & $3$ & $2$ & $2$ & $e_2^{1,0}$ & $2$ & $2$ & $2$ & $\chi^0_1(E_2)$ & $1$ & $1$ & $1$\\
\hline
$e_1^{0,1}$ & $2$ & $2$ & $2$ & $e_2^{0,1}$ & $2$ & $2$ & $2$ & $\chi^0_2(E_2)$ & $1$ & $1$ & $1$\\
\hline
$e_1^{2,0}$ & $3$ & $2$ & $1$ & $e_2^{2,0}$ & $2$ & $2$ & $1$ & $\chi^0_3(E_2)$ & $0$ & $0$ & $0$\\
\hline
$e_1^{1,1}$ & $6$ & $5$ & $5$ & $e_2^{1,1}$ & $4$ & $4$ & $5$ & $\chi^1_1(E_2)$ & $2$ & $2$ & $3$\\
\hline
$e_1^{0,2}$ & $2$ & $2$ & $2$ & $e_2^{0,2}$ & $2$ & $2$ & $2$ & $\chi^1_2(E_2)$ & $2$ & $2$ & $1$\\
\hline
$e_1^{3,0}$ & $1$ & $1$ & $1$ & $e_2^{3,0}$ & $1$ & $1$ & $1$ & $\chi^1_3(E_2)$ & $0$ & $0$ & $0$\\
\hline
$e_1^{2,1}$ & $6$ & $5$ & $4$ & $e_2^{2,1}$ & $4$ & $4$ & $4$ & $\chi^2_1(E_2)$ & $2$ & $2$ & $3$\\
\hline
$e_1^{1,2}$ & $6$ & $5$ & $4$ & $e_2^{1,2}$ & $4$ & $4$ & $4$ & $\chi^2_2(E_2)$ & $2$ & $2$ & $2$\\
\hline
$e_1^{0,3}$ & $1$ & $1$ & $1$ & $e_2^{0,3}$ & $1$ & $1$ & $1$ & $\chi^2_3(E_2)$ & $0$ & $0$ & $0$\\
\hline
$e_1^{3,1}$ & $2$ & $2$ & $2$ & $e_2^{3,1}$ & $2$ & $2$ & $2$ & $\chi^3_1(E_2)$ & $1$ & $1$ & $1$\\
\hline
$e_1^{2,2}$ & $6$ & $5$ & $5$ & $e_2^{2,2}$ & $4$ & $4$ & $5$ & $\chi^3_2(E_2)$ & $1$ & $1$ & $1$\\
\hline
$e_1^{1,3}$ & $3$ & $2$ & $1$ & $e_2^{1,3}$ & $2$ & $2$ & $1$ & $\chi^3_3(E_2)$ & $0$ & $0$ & $0$\\
\hline
$e_1^{3,2}$ & $2$ & $2$ & $2$ & $e_2^{3,2}$ & $2$ & $2$ & $2$ &  &  &  & \\
\hline
$e_1^{2,3}$ & $3$ & $2$ & $2$ & $e_2^{2,3}$ & $2$ & $2$ & $2$ &  &  &  & \\
\hline
\end{tabular}
\end{center}
\end{table}

Since $d_r^{p,q}=0$ for any $(p,q)$ and $r>1$, we have $E_2^{p,q}=E_\infty^{p,q}$. Furthermore, we see from Table $2$ that $\chi^1_1(E_2)$ is lower semi-continuous while $\chi^1_2(E_2)$ is upper semi-continuous.

\bibliographystyle{alpha}
\bibliography{reference}
\end{document}